\documentclass[a4paper]{article}
\usepackage{amssymb,amscd,amsfonts,mathrsfs,amsmath}
\usepackage[all]{xy}
\input amssym.def

\def\double{\mathbb}

\def\cc{{\double C}}
\def\nn{{\double N}}
\def\zz{{\double Z}}

\def\rr{{\double R}}

\newtheorem{theorem}{Theorem}[section]
\newtheorem{lemma}[theorem]{Lemma}
\newtheorem{corollary}[theorem]{Corollary}
\newtheorem{definition}[theorem]{Definition}
\newtheorem{proposition}[theorem]{Proposition}
\newtheorem{remark}[theorem]{Remark}

\newtheorem{example}[theorem]{Example}

\def\si{\sigma}

\def\cinf{C^{\infty}}
\def\cinfc{C^{\infty}_c}

\newcommand{\be}{\begin{equation}}
\newcommand{\ee}{\end{equation}}
\newcommand{\beq}{\begin{eqnarray}}
\newcommand{\eeq}{\end{eqnarray}}
\newcommand{\om}{\omega}
\newcommand{\Om}{\Omega}
\newcommand{\al}{\alpha}

\newcommand{\la}{\lambda}

\newcommand{\Lc}{{\mathscr L}}
\newcommand{\non}{\nonumber}
\newcommand{\eps}{\varepsilon}
\newcommand{\Sc}{{\mathscr S}}

\newcommand{\Ind}{{\mathop{\mathrm{Ind}}}}
\def\ch{\mathrm{ch}}

\def\Td{\mathrm{Td}}

\newcommand{\Tr}{{\mathop{\mathrm{Tr}}}}

\newcommand{\tr}{{\mathop{\mathrm{tr}}}}
\newcommand{\Ac}{{\mathscr A}}

\newcommand{\cqfd}{\hfill\rule{1ex}{1ex}}

\def\Id{\mathrm{Id}}

\def\d{\partial}

\def\im{\mathop{\mathrm{Im}}}

\def\psib{\bar{\psi}}

\def\End{{\mathop{\mathrm{End}}}}

\def\Tc{{\mathscr T}}
\def\Dc{{\mathscr D}}

\def\mod{\ \mathrm{mod}\ }

\def\bint{- \mspace{-20mu} \int}

\def\nablab{\overline{\nabla}}
\def\CL{\mathrm{CL}}
\def\L{\mathrm{L}}
\def\CS{\mathrm{CS}}
\def\PS{\mathrm{PS}}
\def\SPS{\mathrm{SPS}}
\def\i{\mathrm{i}}

\def\Vect{\mathrm{Vect}}

\begin{document}

\begin{center}

{\bf PSEUDODIFFERENTIAL EXTENSION AND TODD CLASS}
\vskip 1cm
{\bf Denis PERROT}
\vskip 0.5cm
Universit\'e de Lyon, Universit\'e Lyon 1,\\
CNRS, UMR 5208 Institut Camille Jordan,\\
43, bd du 11 novembre 1918, 69622 Villeurbanne Cedex, France \\[2mm]
{\tt perrot@math.univ-lyon1.fr}\\[2mm]
\today
\end{center}
\vskip 0.5cm
\begin{abstract}
Let $M$ be a closed manifold. Wodzicki shows that, in the stable range, the cyclic cohomology of the associative algebra of pseudodifferential symbols of order $\leq 0$ is isomorphic to the homology of the cosphere bundle of $M$. In this article we develop a formalism which allows to calculate that, under this isomorphism, the Radul cocycle corresponds to the Poincar\'e dual of the Todd class. As an immediate corollary we obtain a purely algebraic proof of the Atiyah-Singer index theorem for elliptic pseudodifferential operators on closed manifolds.
\end{abstract}

\vskip 0.5cm


\section{Introduction}

Let $M$ be a closed, not necessarily orientable, smooth manifold and denote by $\CL(M)$ the algebra of classical, one-step polyhomogeneous pseudodifferential operators on $M$. The space of smoothing operators $\L^{-\infty}(M)$ is a two-sided ideal in $\CL(M)$, and we call the quotient $\CS(M)=\CL(M)/\L^{-\infty}(M)$ the algebra of \emph{formal symbols} on $M$. The multiplication on $\CS(M)$ is the usual $\star$-product of symbols. One thus gets an extension of associative algebras
\be
0 \to \L^{-\infty}(M) \to \CL(M) \to \CS(M) \to 0\ . \label{ex}
\ee
An ``abstract index problem" then amounts to the computation of the corresponding excision map $HP^{\bullet}(\L^{-\infty}(M)) \to HP^{\bullet+1}(\CS(M))$ in periodic cyclic cohomology \cite{Ni}. In even degree, $HP^0(\L^{-\infty}(M))\cong \cc$ is generated by the usual trace of smoothing operators, whereas in odd degree $HP^1(\L^{-\infty}(M))\cong 0$. Using zeta-function renormalization, one shows (see for instance \cite{P8}) that the image of the trace under the excision map is represented by the following cyclic one-cocycle over $\CS(M)$,
\be
c(a_0,a_1) = \bint a_0  [\log q,a_1] \label{radul}
\ee
for any two formal symbols $a_0,a_1\in\CS(M)$. Here the bar integral denotes the Wodzicki residue \cite{W87}, which is a trace on $\CS(M)$, and $\log q$ is a log-polyhomogeneous symbol associated to a fixed positive elliptic symbol $q\in \CS(M)$ of order one. Notice that the bilinear functional $c$ was originally introduced by Radul in the context of Lie algebra cohomology \cite{Ra}. A direct computation shows that $c$ is in fact a cyclic one-cocycle over $\CS(M)$, and that its cyclic cohomology class does not depend on the choice of $q$. Hence the class $[c]\in HP^1(\CS(M))$ is completely canonical, in the sense that it only depends on $M$. On the other hand the cyclic cohomology of $\CS(M)$ is known \cite{W88}, and corresponds to the ordinary homology (with complex coefficients) of a certain manifold. A natural question therefore is to identify the class $[c]$. In the present paper we give the answer for its image in the periodic cyclic cohomology of the subalgebra $\CS^0(M)\subset \CS(M)$, the formal symbols of order $\leq 0$. The result is stated as follows. The leading symbol map gives rise to an algebra homomorphism $\la: \CS^0(M)\to \cinf(S^*M)$, where $S^*M$ is the cosphere bundle of $M$. This allows to pullback any homology class of $S^*M$ to the periodic cyclic cohomology of the symbol algebra:
\be
\la^*: H_{\bullet}(S^*M,\cc) \to HP^{\bullet}(\CS^0(M))\ .
\ee
Wodzicki shows that $\la^*$ is an \emph{isomorphism}, provided that the natural locally convex topology of $\CS^0(M)$ is taken into account \cite{W88}. Our main result is the following theorem (\ref{trad}), which holds in the algebraic setting or the locally convex setting regardless to Wodzicki's isomorphism.

\begin{theorem}\label{tt}
Let $M$ be a closed manifold. The periodic cyclic cohomology class of $[c]\in HP^1(\CS^0(M))$ is
\be
[c] =  \la^*\big([S^*M]\cap\, \pi^*\Td(T_{\cc}M)\big)\ , \label{equ}
\ee
where $\Td(T_{\cc}M) \in H^{\bullet}(M,\cc)$ is the Todd class of the complexified tangent bundle,  and $\pi:S^*M\to M$ is the cosphere bundle endowed with its canonical orientation and fundamental class $[S^*M]\in H_{\bullet}(S^*M)$. 
\end{theorem}
We give a purely algebraic proof of this theorem. The central idea is to introduce the $\zz_2$-graded algebra $\CL(M,E)$ of pseudodifferential operators acting on differential forms, that is, on the sections of the exterior bundle $E=\Lambda T^*M$, and view the corresponding algebra of formal symbols $\CS(M,E)$ as a bimodule over itself. Using this bimodule structure we develop a formalism of abstract Dirac operators. This leads to the construction of cyclic cocycles for the subalgebra $\CS^0(M)\subset \CS(M,E)$. These cocycles are given by algebraic analogues of the JLO formula \cite{JLO}, and are all cohomologous in $HP^{\bullet}(\CS^0(M))$. By choosing genuine Dirac operators we obtain both sides of equality (\ref{equ}). Let us mention that the JLO formula in the right-hand-side provides a representative of the Todd class as a closed differential form over $M$
\be
\Td(\i R/2\pi) = \det\left(\frac{\i R/2\pi}{e^{\i R/2\pi}-1}\right)
\ee
where $R$ is the curvature two-form of an affine torsion-free connection on $M$. Hence our method gives an ``explicit formula" for the class $[c]$. In the same way, we also prove that the cyclic cohomology class of the Wodzicki residue vanishes in $HP^0(\CS^0(M))$.\\

As an immediate corollary of Theorem \ref{tt} we obtain the Atiyah-Singer index formula for elliptic pseudodifferential operators \cite{AS}. If $Q$ is an elliptic operator acting on the sections of a (trivially graded) vector bundle over $M$, its leading symbol is an invertible matrix $g$ with entries in the commutative algebra $\cinf(S^*M)$, hence it defines a class in the algebraic $K$-theory $K_1(\cinf(S^*M))$. Its Chern character in $H^{\bullet}(S^*M,\cc)$ is represented by the closed differential form of odd degree
\be
\ch(g) = \sum_{k\geq 0}  \frac{k!}{(2k+1)!} \, \tr\Big(\frac{(g^{-1}dg)^{2k+1}}{(2\pi\i)^{k+1}}\Big)\ .
\ee
\begin{corollary}[Index theorem] 
Let $Q$ be an elliptic pseudodifferential operator of order $\leq 0$ acting on the sections of a trivially graded vector bundle over $M$, with leading symbol class $[g]\in K_1(\cinf(S^*M))$. Then the Fredholm index of $Q$ is the integer
\be
\Ind(Q) = \langle [S^*M], \pi^*\Td(T_{\cc}M)\cup \ch([g])\rangle\ .
\ee
\end{corollary}
This is a direct consequence of the fact that the class $[c]\in HP^1(\CS^0(M))$ of the residue cocycle is the image of the operator trace $\Tr:\L^{-\infty}(M)\to\cc$ under the excision map of the fundamental extension 
\be
0 \to \L^{-\infty}(M) \to \CL^0(M) \to \CS^0(M) \to 0 \ .
\ee
In fact (\ref{equ}) and the index formula are equivalent. Hence our method gives a new algebraic proof of the index theorem. This should however not be confused with what is usually called an \emph{algebraic index theorem} (\cite{NT}). The latter calculates the cyclic cohomology class of the canonical trace on a (formal) deformation quantization of the algebra of smooth functions on a symplectic manifold, and relates it to the Todd class. In the special case of the symplectic manifold $T^*M$, one may take the algebra of smoothing operators $\L^{-\infty}(M)$ as a deformation quantization of the commutative algebra of functions over $T^*M$ and obtain in this way the usual index theorem. This is \emph{not} what we are doing here. In fact our approach is in some sense opposite, because instead of working with the operator ideal $\L^{-\infty}(M)\subset \CL^0(M)$ we directly deal with the quotient algebra of formal symbols $\CS^0(M)$. As a consequence, we drop the delicate analytic issues inherent to the highly non-local algebra $\L^{-\infty}(M)$ and its operator trace, and entirely transfer the index problem on the algebra $\CS^0(M)$ endowed with the residue cocycle (\ref{radul}). The computation is purely local because only a finite number of terms in the asymptotic expansion of symbols contribute to the index, which relates our approach to the Connes-Moscovici residue index formula \cite{CM95}. For this reason our formalism is well-adapted (and in fact motivated by) the study of more general index problems appearing in non-commutative geometry \cite{C94}, for which a genuine extension of algebras and the corresponding residue cocycle are available. This includes higher equivariant index theorems for non-isometric actions of non-compact groups, higher index theorems on Lie groupoids, and so on. These ideas will be developped elsewhere.\\

Here is a brief description of the paper. In section \ref{spo} we recall basic things about pseudodifferential operators. In section \ref{sbfs} we look at $\CS(M,E)$ as a bimodule over itself and introduce the relevant spaces of operators acting on it. In section \ref{sct} a canonical trace is defined by means of the Wodzicki residue. Section \ref{sdo} introduces generalized Dirac operators acting on $\CS(M,E)$. Theorem \ref{tt} is proved in section \ref{sjlo} by means of the algebraic JLO formula, and the index theorem is deduced in section \ref{sit}. \\
All manifolds are supposed to be Hausdorff, paracompact, smooth and without boundary.

\section{Pseudodifferential operators}\label{spo}

Let $M$ be a $n$-dimensional manifold. We denote by $\cinf(M)$ (resp. $\cinfc(M)$) the space of smooth complex-valued (resp. compactly supported) functions over $M$. A linear map $A:\cinfc(M)\to\cinf(M)$ is a pseudodifferential operator of order $m\in \rr$ if for every coordinate chart $(x^1,\ldots,x^n)$ over an open subset $U\subset M$, there exists a smooth function $a\in \cinf(U\times \rr^n)$ such that
\be
(A\cdot f)(x) = \frac{1}{(2\pi)^n}\int_{U\times\rr^n} e^{\i p\cdot(x-y)}a(x,p) \, f(y) \, dy\, dp
\ee
for any $f\in \cinfc(U)$. We use the notation $\i=\sqrt{-1}$. For any multi-indices $\al=(\al_1,\ldots,\al_n)$ and $\beta=(\beta_1,\ldots,\beta_n)$, the symbol $a$ has to satisfy the estimate
\be
|\d_x^{\al}\d_p^{\beta} a(x,p)| \leq C_{\al,\beta}(1+\|p\|)^{m-|\beta|}
\ee
for some constant $C_{\al,\beta}$, where $|\beta|=\beta_1+\ldots+\beta_n$, $\d_x=\frac{\d}{\d x}$ and $\d_p=\frac{\d}{\d p}$ are the partial derivatives with respect to the variables $x=(x^1,\ldots,x^n)$ and $p=(p_1,\ldots,p_n)$, and $\|p\|$ is the euclidian norm of $p\in \rr^n$. Note that $(x,p)$ is the canonical coordinate system on the cotangent bundle $T^*U\cong U\times \rr^n$. In addition, $A$ is a \emph{classical} (one-step polyhomogeneous) pseudodifferential operator of order $m$ if its symbol in any coordinate chart has an asymptotic expansion as $\|p\|\to\infty$ of the form
\be
a(x,p) \sim \sum_{j=0}^{\infty} a_{m-j}(x,p) \label{expan}
\ee
where the functions $a_{m-j}\in \cinf(U\times\rr^n)$ are homogeneous of degree $m-j$ with respect to the variable $p$. For any $m\in \rr$, we denote by $\CL^m(M)$ the space of all classical pseudodifferential operators of order $\leq m$. One has $\CL^m(M)\subset \CL^{m'}(M)$ whenever $m\leq m'$. Define as usual the space of all classical pseudodifferential operators and the space of smoothing operators, respectively
\be
\CL(M) = \bigcup_{m\in\rr} \CL^m(M)\ ,\qquad \L^{-\infty} = \bigcap_{m\in\rr} \CL^m(M) \ .
\ee
Two operators in $\CL(M)$ are equal modulo smoothing operators if and only if their asymptotic expansions (\ref{expan}) agree in all coordinate charts. The space of \emph{formal classical symbols} $\CS(M)$ is defined via the exact sequence
\be
0 \to \L^{-\infty}(M) \to \CL(M) \to \CS(M) \to 0 \label{ext}
\ee
Thus, a formal symbol of order $m$ corresponds to a formal series as the right-hand-side of (\ref{expan}) in any local chart, which fulfills complicated gluing formulas under coordinate-change. $\CS(M)$ is of course the union, for all $m\in\rr$, of the subspaces $\CS^m(M)$ of formal symbols of order $\leq m$. Recall that $\CS^m(M)$ is \emph{complete}, in the sense that any formal series of homogeneous functions $a_{m-j}$ is the formal symbol of some pseudodifferential operator. We denote by $\PS(M)\subset \CS(M)$ the subalgebra of formal symbols which are polynomial with respect to the variable $p$ in any chart. $\PS(M)$ is isomorphic to the space of differential operators on $M$. \\
The composition of pseudodifferential operators is not always defined, unless these operators are properly supported. This happens in particluar when $M$ is compact. In that case, $\CL(M)$ becomes a filtered associative algebra, i.e. $\CL^m(M)\cdot\CL^{m'}(M)\subset \CL^{m+m'}(M)$, and $\L^{-\infty}(M)$ is a two-sided ideal. Hence (\ref{ext}) is actually an exact sequence of associative algebras. The product of two formal symbols $a,b\in\CS(M)$ in a local chart is the $\star$-product
\be
(ab)(x,p) = \sum_{|\al|=0}^{\infty} \frac{(-\i)^{|\al|}}{\al !}\, \d_p^{\al}a(x,p) \, \d_x^{\al}b(x,p) \label{star}
\ee
where the sum runs over all multi-indices $\al=(\al_1,\ldots,\al_n)$, and $\al! = \al_1!\ldots\al_n!$. Notice that, in contrast with $\CL(M)$, the product in $\CS(M)$ is defined without any condition on the support (compact or not) of the symbols. \\

If $E$ is a (possibly $\zz_2$-graded) complex vector bundle over $M$, the algebra of classical pseudodifferential operators $\CL(M,E)$ acting on the smooth sections of $E$ is defined analogously. The only difference is that over a local chart which also trivialises $E$, the symbol becomes a function of $(x,p)$ with values in the matrix algebra $M_k(\cc)$ where $k$ is the rank of $E$. One has the exact sequence
\be
0 \to \L^{-\infty}(M,E) \to \CL(M,E) \to \CS(M,E) \to 0\ .
\ee
$\PS(M,E)\subset \CS(M,E)$ denotes the algebra of polynomial symbols with respect to $p$. It is isomorphic to the algebra of differential operators acting on the smooth sections of $E$. In the sequel we will essentially focus on the $\zz_2$-graded bundle $E=\Lambda T_{\cc}^*M$, the exterior algebra of the complexified cotangent bundle of $M$. The smooth sections of $E$ are the complex differential forms over $M$. Consider the (real) vector bundle $TM\oplus T^*M$ endowed with its canonical inner product. Then $E$ is a spinor representation of the complexified Clifford algebra bundle $C(TM\oplus T^*M)$. In other words, the endomorphism bundle $\End(E)$ is canonically isomorphic to $C(TM\oplus T^*M)$. We use this identification in order to find a set of generators for the algebra $\PS(M,E)$ in a local coordinate system $(x^1,\ldots,x^n)$ over an open $U\subset M$. For each $i$ we view $x^i$ as the multiplication operator of a differential form by the function $x^i$, and $\i p_i$ as the Lie derivative of a differential form with respect to the vector field $\frac{\d}{\d x^i}$. For all indices $i,j$ they fulfill the usual Canonical Commutation Relations
\be
[x^i,x^j]=0\ ,\qquad [x^i,p_j]=\i \delta^i_j\ ,\qquad [p_i,p_j]=0
\ee
($\i =\sqrt{-1}$), and generate the even part of the algebra of differential operators $\PS(U,E)$. The odd generators are defined by the operators 
\be
\psi^i = \mu(dx^i)\ ,\qquad \psib_i = \iota(\d_{x^i})\ ,
\ee
where $\mu$ is exterior multiplication by a differential form (on the left) and $\iota$ is interior multiplication by a vector field (on the left). These are odd sections of the endomorphism bundle $\End(E)$ over $U$, and their \emph{graded} commutators (hence anticommutators) fulfill the Clifford relations (Canonical Anticommutation Relations)
\be
[\psi^i,\psi^j] =0\ ,\qquad [\psi^i,\psib_j] =\delta^i_j\ ,\qquad [\psib_i,\psib_j] =0 \label{clif}
\ee
while the commutators between $x,p$ on one hand and $\psi,\psib$ on the other hand all vanish. The odd operators $\psi,\psib$ generate a basis of sections for $\End(E)$ over $U$. Hence a differential operator $a\in \PS(M,E)$ is represented over $U$ as a function $a(x,p,\psi,\psib)$ which depends polynomially on the even variable $p$. Since the odd variables generate a finite-dimensional algebra, $a$ is also a polynomial with respect to $\psi,\psib$. In the same way, any symbol $a\in \CS(M,E)$ of order $m$ is locally represented as a formal series, over $j\in\nn$, of functions $a_{m-j}(x,p,\psi,\psib)$ which are homogeneous of degree $m-j$ with respect to $p$ and polynomial with respect to the odd variables $\psi,\psib$. \\
Let us end this paragraph with the effect of a coordinate change (or local diffeomorphism) $\gamma$ on the generators $(x,p,\psi,\psib)$ of $\CS(U,E)$. If one puts $\gamma(x^i)=y^i$ for all $i$, then 
\be
\gamma(\psi^i)= \mu(dy^i)=\mu\Big(\frac{\d y^i}{\d x^j}\, dx^j\Big) = \frac{\d y^i}{\d x^j}\, \psi^j
\ee
where we use Einstein's convention of summation over repeated indices. In the same way
\be
\gamma(\psib_i) = \iota(\d_{y^i}) = \iota\Big(\frac{\d x^j}{\d y^i}\, \d_{x^j}\Big) = \frac{\d x^j}{\d y^i} \, \psib_j\ .
\ee
Finally, the identification of $\i p_i$ with the Lie derivative $\iota(\d_{x^i})\circ d + d \circ\iota(\d_{x^i}) = \psib_i \circ d + d\circ \psib_i$ yields
$$
\gamma(p_i) = -\i (\gamma(\psib_i)\circ d + d\circ \gamma(\psib_i)) = \frac{\d x^j}{\d y^i}\, p_j - \i \, d\Big(\frac{\d x^j}{\d y^i}\Big)\psib_j\ .
$$
Since the exterior derivative is $d=\psi^k\d_{x^k}$, and $\psi^k$ commutes with functions of $x$, one has
\be
\gamma(p_i) = \frac{\d x^j}{\d y^i}\, p_j - \i \, \frac{\d}{\d x^k}\Big(\frac{\d x^j}{\d y^i}\Big)\psi^k\psib_j\ .
\ee

\section{The bimodule of formal symbols}\label{sbfs}

Let $M$ be an $n$-dimensional manifold and consider the $\zz_2$-graded algebra of formal symbols $\CS(M,E)$ with $E=\Lambda T^*_{\cc}M$. We view $\CS(M,E)$ as a left $\CS(M,E)$-module and right $\PS(M,E)$-module: the left action of a formal symbol $a\in\CS(M,E)$ and the right action of a polynomial symbol $b\in \PS(M,E)$ on a vector $\xi \in \CS(M,E)$ read
\be
a_L\cdot \xi = a \xi\ ,\qquad b_R\cdot \xi = \pm\, \xi b\ ,
\ee
where the sign $\pm$ depends on the respective parities of $b$ and $\xi$: it is $-$ if both $b$ and $\xi$ are odd, $+$ otherwise. The left action of $a$ induces a representation of $\CS(M,E)$ in the algebra of linear endomorphisms $\End(\CS(M,E))$. The right action of $b$ induces a representation of the opposite algebra $\PS(M,E)^{\mathrm{op}}$ in $\End(\CS(M,E))$. The left and right actions commute in the graded sense, whence an algebra homomorphism from the (graded) tensor product $\CS(M,E)\otimes\PS(M,E)^{\mathrm{op}}$ to $\End(\CS(M,E))$. This homomorphism is not injective. Its range defines a $\zz_2$-graded algebra 
\be
\Lc(M) = \im \big( \CS(M,E)\otimes \PS(M,E)^{\mathrm{op}}\to \End(\CS(M,E)) \big)\ .
\ee
Thus $\Lc(M)$ is linearly generated by products $a_Lb_R$ with $a\in\CS(M,E)$ and $b\in \PS(M,E)$. Let $(x^1,\ldots,x^n)$ be a local coordinate system over an open subset $U\subset M$. The function $x^i$ is a symbol (of order zero) in $\PS(U,E)$, so that $x^i_L$ and $x^i_R$ are elements of $\Lc(U)$. Moreover for any $\xi\in \CS(U,E)$ one has
\be
(x^i_L-x^i_R)\cdot\xi = [x^i,\xi] = \i \frac{\d\xi}{\d p_i}\ ,
\ee
whence $x^i_L-x^i_R= \i\frac{\d}{\d p_i}$. In the same way the conjugate coordinate $p_i$ is a symbol (of order one) in $\PS(U,E)$, so $p_{iL}$ and $p_{iR}$ are elements of $\Lc(U)$, and for any $\xi\in \CS(U,E)$,
\be
(p_{iL}-p_{iR})\cdot\xi = [p_i,\xi] = -\i \frac{\d\xi}{\d x^i}\ ,
\ee
whence $p_{iL}-p_{iR}= -\i\frac{\d}{\d x^i}$. The situation is analogous for the odd coordinates $\psi^i$ and $\psib_i$, and one finds that $\psi^i_L-\psi^i_R$ is the partial derivative with respect to $\psib_i$, while $\psib_{iL}-\psib_{iR}$ is the partial derivative with respect to $\psi^i$. If $b\in\PS^k(M,E)$ is a differential operator of order $k\in\nn$, we can write, locally over $U$
$$
b(x,p,\psi,\psib) = \sum_{|\al|=0}^{k} b_{\al}(x,\psi,\psib)\, p^{\al} = \sum_{|\al|=0}^{k} \sum_{|\eta|=0}^{n} \sum_{|\theta|=0}^{n} b_{\al,\eta,\theta}(x)\, p^{\al}\, \psi^{\eta} \, \psib^{\theta} \ ,
$$
where $b_{\al,\eta,\theta}$ are functions of the only variable $x$, and $\al,\eta,\theta$ are multi-indices. Using formula (\ref{star}) for the star-product, one finds  
$$
(b_{\al,\eta,\theta})_R \cdot \xi =  \sum_{|\beta|=0}^{\infty} \frac{(-\i)^{|\beta|}}{\beta !}\, (\d_x^{\beta}b_{\al,\eta,\theta})_L \cdot \d_p^{\beta} \xi  
$$
for any $\xi\in\CS(M,E)$. Since left and right actions commute, the operator $b_R$ reads
$$
b_R = \sum_{|\al|=0}^{k} \sum_{|\beta|=0}^{\infty} \sum_{|\eta|=0}^{n} \sum_{|\theta|=0}^{n} \frac{(-\i)^{|\beta|}}{\beta !}\,  (\d_x^{\beta}b_{\al,\eta,\theta})_L \, (\psi^{\eta}\psib^{\theta})_R   \, p_R^{\al} \, \d_p^{\beta} 
$$
Using the identity $p_R=p_L+\i\d_x$, one concludes that a generic element $a_Lb_R \in\Lc(M)$ can be expressed, locally over a subset $U\subset M$, as a series
\be
a_Lb_R = \sum_{|\al|=0}^{k} \sum_{|\beta|=0}^{\infty} \sum_{|\eta|=0}^{n} \sum_{|\theta|=0}^{n} (s_{\al,\beta,\eta,\theta})_L \, (\psi^{\eta}\psib^{\theta})_R \, \d_x^{\al} \, \d_p^{\beta} \ , \label{generic}
\ee
for some coefficients $s_{\al,\beta,\eta,\theta} \in \CS(U,E)$ and finite $k\in\nn$. Notice however that the converse is not true: a series (\ref{generic}) with arbitrary coefficients $s_{\al,\beta,\eta,\theta}$ does not necessarily come from an element of $\Lc(M)$. \\ 

Now let $\Sc(M)=\Lc(M)[[\eps]]$ be the $\zz_2$-graded algebra of \emph{formal power series} in the indeterminate $\eps$, with coefficients in $\Lc(M)$. The generator $\eps$ has trivial grading. An element of $\Sc(M)$ is therefore an infinite sum $s=\sum_{k=0}^{\infty} s_k\eps^k$ where each coefficient $s_k$ is given by a series of the form (\ref{generic}) in any local chart. We can view $\Sc(M)$ as an algebra of linear operators acting on the space of formal series $\CS(M,E)[[\eps]]$. This algebra is filtered by the subalgebras $\Sc_k(M)=\Sc(M)\eps^k$, $\forall k\in\nn$. For each $m\in\rr$, we define a subspace $\Dc^m(M)\subset \Sc(M)$ as follows. An element $s=\sum s_k\eps^k$ is in $\Dc^m(M)$ if and only if in any local chart over $U\subset M$,
\be
s_k = \sum_{|\al|=0}^{k} \sum_{|\beta|=0}^{\infty} \sum_{|\eta|=0}^{n} \sum_{|\theta|=0}^{n} (s^m_{k,\al,\beta,\eta,\theta})_L \, (\psi^{\eta}\psib^{\theta})_R \, \d_x^{\al} \, \d_p^{\beta} 
\ee
where $s^m_{k,\al,\beta,\eta,\theta} \in \CS(U,E)$ is a symbol of order $\leq m+(k+|\beta|-3|\al|)/2$. Moreover we set $\Dc_k^m(M)=\Dc^m(M)\cap \Sc_k(M)$ for all $k\in\nn$. Hence the subscript $k$ counts the minimal power of $\eps$ appearing in a formal series. Observe that in local coordinates, the partial derivative $\d_x$ always appears with at least one power of $\eps$. Here are some examples: $1\in \Dc^0_0(M)$, $\CS^m(M,E)_L\subset \Dc_0^m(M)$, $\eps\in\Dc^{-1/2}_1(M)$, $\eps\d_x\in \Dc_1^1(U)$, $\d_p\in\Dc_0^{-1/2}(U)$, and   $\eps\d_x\d_p\in \Dc_1^{1/2}(U)$. One obviously has $\Dc^m(M)\subset\Dc^{m'}(M)$ whenever $m\leq m'$, and we set $\Dc(M)=\bigcup_{m\in\rr}\Dc^m(M)$. The following lemma shows that $\Dc(M)$ is a subalgebra of $\Sc(M)$.

\begin{lemma}
The inclusion $\Dc_k^m(M)\Dc_{k'}^{m'}(M)\subset \Dc_{k+k'}^{m+m'}(M)$ holds in all degrees $m,m'\in \rr$ and $k,k'\in\nn$. Hence $\Dc(M)$ is a unital, $\zz_2$-graded, bi-filtered algebra. 
\end{lemma}
{\it Proof:} Since $\psi_R$ and $\psib_R$ play no role in the filtration degrees, it suffices to show that, in a local coordinate system over $U$, the composition $s\circ s'$ of two operators
\beq
s &=& \sum_{k=0}^{\infty} \sum_{|\al|=0}^{k} \sum_{|\beta|=0}^{\infty} \eps^k(s^m_{k,\al,\beta})_L \, \d_x^{\al}  \d_p^{\beta} \ \in \Dc^m(U)\ , \non\\
s'&=& \sum_{k'=0}^{\infty} \sum_{|\al'|=0}^{k'} \sum_{|\beta'|=0}^{\infty} \eps^{k'}(s^{m'}_{k',\al',\beta'})_L \, \d_x^{\al'}  \d_p^{\beta'} \ \in \Dc^{m'}(U)\non
\eeq
is in $\Dc^{m+m'}(U)$. Note that the commutator $[\d_p,\ ]$ on a symbol decreases the order by one, whereas the commutator $[\d_x, \ ]$ leaves the order unaffected. Hence we can write the composition $\d_x^{\al} \d_p^{\beta} \circ (s^{m'}_{k',\al',\beta'})_L$ as a sum
$$
\d_x^{\al} \d_p^{\beta} \circ (s^{m'}_{k',\al',\beta'})_L = \sum_{|\gamma|=0}^{|\al|} \sum_{|\delta|=0}^{|\beta|} (t^{m',\al,\beta}_{k',\al',\beta',\gamma,\delta})_L \d_x^{\gamma} \d_p^{\delta}
$$
where $t^{m',\al,\beta}_{k',\al',\beta',\gamma,\delta}$ is a symbol of order $\leq m'-|\beta|+|\delta|+(k'+|\beta'|-3|\al'|)/2 $. Then
$$
s\circ s' = \sum_{k,k',|\beta|,|\beta'| \geq 0} \sum_{|\al|=0}^k \sum_{|\al'|=0}^{k'} \sum_{|\gamma|=0}^{|\al|} \sum_{|\delta|=0}^{|\beta|} \eps^{k+k'} \big(s^m_{k,\al,\beta} t^{m',\al,\beta}_{k',\al',\beta',\gamma,\delta}\big)_L \d_x^{\gamma+\al'} \d_p^{\delta+\beta'}
$$
The symbol $s^m_{k,\al,\beta} t^{m',\al,\beta}_{k',\al',\beta',\gamma,\delta}$ has order $\leq m+m'-|\beta|+|\delta|+\frac{1}{2}(k+k'+|\beta|+|\beta'|-3|\al|-3|\al'|) = m+m' + \frac{1}{2}(k+k'+|\delta+\beta'|-3|\gamma+\al'|) -\frac{3}{2}(|\al|-|\gamma|) -\frac{1}{2}(|\beta|-|\delta|) $. For fixed indices $k,k',\al',\beta',\gamma,\delta$ this order is a strictly decreasing function of $|\al|$ and $|\beta|$. Moreover $\frac{3}{2}(|\al|-|\gamma|)\geq 0$ and $\frac{1}{2}(|\beta|-|\delta|)\geq 0$. Hence by completeness of the space of symbols, the series
$$
u^{m+m'}_{k,k',\al',\beta',\gamma,\delta} = \sum_{|\al|=|\gamma|}^{k} \sum_{|\beta|=0}^{\infty} s^m_{k,\al,\beta} t^{m',\al,\beta}_{k',\al',\beta',\gamma,\delta}
$$
converges to a symbol of order $\leq m+m' + \frac{1}{2}(k+k'+|\delta+\beta'|-3|\gamma+\al'|)$. It follows that
$$
s\circ s' = \sum_{k,k' \geq 0}  \sum_{|\al'|=0}^{k'} \sum_{|\gamma|=0}^{k} \sum_{|\beta|=0}^{\infty} \sum_{|\delta|=0}^{|\beta|} \eps^{k+k'} \big(u^{m+m'}_{k,k',\al',\beta',\gamma,\delta}\big)_L \d_x^{\gamma+\al'} \d_p^{\delta+\beta'}
$$
is indeed an element of $\Dc^{m+m'}(U)$. This shows the inclusion $\Dc^m(M)\Dc^{m'}(M)\subset \Dc^{m+m'}(M)$. Furthermore $\Sc_k(M) \Sc_{k'}(M)\subset \Sc_{k+k'}(M)$ is obvious, one concludes that $\Dc_k^m(M)  \Dc_{k'}^{m'}(M) \subset \Dc_{k+k'}^{m+m'}(M)$. \cqfd\\

\begin{definition}\label{dlap}
An operator $\Delta\in\Dc^{1/2}_1(M)$ of even parity is called a \emph{generalized Laplacian} if in any coordinate system over an open set $U\subset M$ it reads
\be
\Delta \equiv \i\eps \frac{\d}{\d x^i}\frac{\d}{\d p_i} \mod \Dc^0_1(U)
\ee
(summation over repeated indices).
\end{definition}
 
\begin{lemma}
A generalized Laplacian exists over any manifold $M$.
\end{lemma}
{\it Proof:} It is actually a consequence of the existence of Dirac operators (section \ref{sdo}) but we can give a direct proof by looking at the behaviour of the canonical ``flat" Laplacian $\i\eps \frac{\d}{\d x^i}\frac{\d}{\d p_i}$ under a coordinate change $x^i\mapsto \gamma(x^i)=y^i$ over $U$. One has $\frac{\d}{\d x^i} = \i (p_{iL}- p_{iR}$) hence
$$
\gamma\Big(\frac{\d}{\d x^i}\Big) = \i \gamma(p_i)_L -\i \gamma(p_i)_R\ .
$$
Recall that $\gamma(p_i) = \frac{\d x^j}{\d y^i}\, p_j - \i \, \frac{\d}{\d x^k}\frac{\d x^j}{\d y^i}\,\psi^k\psib_j$. The operators $(\frac{\d}{\d x^k}\frac{\d x^j}{\d y^i}\,\psi^k\psib_j)_L$ and $(\frac{\d}{\d x^k}\frac{\d x^j}{\d y^i}\,\psi^k\psib_j)_R$ belong to $\Dc^0(U)$, so that
$$
\gamma\Big(\i\eps\frac{\d}{\d x^i}\Big) \equiv - \eps\Big(\frac{\d x^j}{\d y^i}p_j\Big)_L + \eps\Big(\frac{\d x^j}{\d y^i}p_j\Big)_R \mod \Dc^0_1(U)\ .
$$
Then we use the expansion 
$$
\eps\Big(\frac{\d x^j}{\d y^i}p_j\Big)_R = \eps \sum_{|\beta|=0}^{\infty} \frac{(-\i)^{|\beta|}}{\beta !} \Big(\d_x^{\beta}\frac{\d x^j}{\d y^i}\Big)_L(p_j)_R \d_p^{\beta}\ .
$$
Since $\eps p_R \in \Dc^1_1(U)$ and $\d_p\in \Dc^{-1/2}_0(U)$, the terms in the right-hand side belong to $\Dc^{1/2}_1(U)$ whenever $|\beta|\geq 1$. Thus we only retain the principal term $|\beta|=0$:
$$
\gamma\Big(\i\eps \frac{\d}{\d x^i}\Big) \equiv \eps \Big(\frac{\d x^j}{\d y^i}\Big)_{\!\!L} \big(-p_{jL} + p_{jR} \big) \equiv \i\eps \Big(\frac{\d x^j}{\d y^i}\Big)_{\!\!L} \frac{\d}{\d x^j} \mod \Dc^{1/2}_1(U) \ .
$$
We proceed in the same way with $\frac{\d}{\d p_i} = -\i x^i_L + \i x^i_R$:
$$
\gamma\Big(\frac{\d}{\d p_i}\Big) = -\i \gamma(x^i)_L + \i \gamma(x^i)_R = -\i y^i_L + \i y^i_R\ .
$$
We use the expansion 
$$
y^i_R = \sum_{|\beta|=0}^{\infty} \frac{(-\i)^{|\beta|}}{\beta !}  (\d_x^{\beta}y^i)_L \d_p^{\beta}\ .
$$
Since $\d_p^{\beta}\in \Dc^{-|\beta|/2}_0(U)$, we only retain the principal terms $|\beta|=1$ in the following sum:
$$
\gamma\Big(\frac{\d}{\d p_i}\Big) = \i \sum_{|\beta|=1}^{\infty} \frac{(-\i)^{|\beta|}}{\beta !}  (\d_x^{\beta}y^i)_L \d_p^{\beta} \equiv \Big(\frac{\d y^i}{\d x^j}\Big)_{\!\!L} \frac{\d}{\d p_j} \mod \Dc^{-1}_0(U) \ .
$$
Finally we can write
\beq
\gamma\Big(\frac{\d}{\d p_i}\Big)\gamma\Big(\i\eps \frac{\d}{\d x^i}\Big) &\equiv& \Big(\Big(\frac{\d y^i}{\d x^j}\Big)_{\!\!L} \frac{\d}{\d p_j} \mod \Dc^{-1}_0\Big) \Big(\i\eps\Big(\frac{\d x^j}{\d y^i}\Big)_{\!\!L} \frac{\d}{\d x^j} \mod \Dc^{1/2}_1\Big) \non\\
&\equiv&  \i\eps\frac{\d}{\d p_j}\frac{\d}{\d x^j}\ \mod \big(\Dc^{-1}_0 \Dc^{1}_1 + \Dc^{-1/2}_0 \Dc^{1/2}_1 \big) \non
\eeq
This shows that the operator $\i\eps \frac{\d}{\d x^i}\frac{\d}{\d p_i}$ is invariant modulo $\Dc^0_1(U)$ under coordinate change. Now let $(c_I)$ be a partition of unity associated to an atlas $(U_{I},x_{I})$ on $M$. Denoting by $\Delta_I$ the canonical flat Laplacian in local coordinates $x_I$, the sum
$$
\Delta = \sum_I (c_I)_L\Delta_I
$$
globally defines a generalized Laplacian on $M$. \cqfd\\

\noindent Observe that a generalized Laplacian $\Delta$ carries at least one power of $\eps$, hence any formal power series of $\Delta$ is a well-defined element of $\Sc(M)$. For example, for any parameter $t\in\rr$ the exponential
\be
\exp(t\Delta) = \sum_{k=0}^{\infty} \frac{t^k}{k!}\, \Delta^k
\ee
is an invertible element of $\Sc(M)$, with inverse $\exp(-t\Delta)$. However, these elements do not belong to $\Dc(M)$. We define an automorphism $\si^t_{\Delta}$ of the algebra $\Sc(M)$ as follows:
\be
\si^t_{\Delta}(s) = \exp(t\Delta)\, s \, \exp(-t\Delta) \qquad \forall s\in\Sc(M) \ .
\ee
Clearly $\si^t_{\Delta}\circ \si^{t'}_{\Delta}=\si^{t+t'}_{\Delta}$ so the map $t\mapsto \si^t_{\Delta}$ defines a one-parameter group of automorphisms.
\begin{lemma} \label{lsi}
For any generalized Laplacian $\Delta$, the automorphism group $\si_{\Delta}$ preserves the subalgebra $\Dc(M)$. More precisely one has $[\Delta,\Dc^m_k(M)]\subset \Dc^m_{k+1}(M)$ and $\si^t_{\Delta}(\Dc_k^m(M))= \Dc_k^m(M)$ for all $m\in\rr$, $k\in\nn$ and $t\in\rr$.
\end{lemma}
{\it Proof:} In local coordinates $\Delta=\i\eps\frac{\d}{\d x^i}\frac{\d}{\d p_i}+r$ with $r \in\Dc^0_1$. Hence for any $s\in\Dc^m_k$
$$
[\Delta,s] = \i\eps (\d_{x}s\, \d_{p} + \d_{p}s\, \d_{x}+\d_{x}\d_{p}s ) + [r,s]\ .
$$
One has $\eps\d_xs\d_p\in \Dc^{m-1}_{k+1}$, $\eps\d_ps\d_x\in \Dc^m_{k+1}$, $\eps\d_x\d_ps\in \Dc^{m-3/2}_{k+1}$, $r s\in \Dc^m_{k+1}$ and $s r \in \Dc^m_{k+1}$. Hence $[\Delta,\Dc^m_k(M)]\subset \Dc^m_{k+1}(M)$ as claimed.\\
Next we show $\exp(t\Delta)\Dc_k^m(M) \exp(-t\Delta) \subset \Dc_k^m(M)$ for all $m,k$. Replacing $t$ by $-t$ then gives the inverse inclusion. For $s\in\Dc^m_k(M)$ consider the identity
$$
\exp(t\Delta)\, s \, \exp(-t\Delta) = \sum_{l=0}^{\infty} \frac{t^l}{l!}\,  s^{(l)}
$$
where $s^{(l)}$ denotes the $l$-th power of the derivation $[\Delta,\ ]$ on $s$. By induction one has $s^{(l)}\in \Dc^m_{k+l}(M)$ for any $l\geq 0$ so that the infinite sum over $l$ gives a well-defined element of $\Dc^m_k(M)$. \cqfd\\
\begin{lemma}
Let $\Delta+s$ be a perturbation of a generalized Laplacian $\Delta$, with $s\in \Dc^0_1(M)$. Then the Duhamel formula holds in $\Sc(M)$:
\be
\exp(\Delta+s) = \sum_{k=0}^{\infty} \int_{\Delta_k} \exp(t_0\Delta)\, s \, \exp(t_1\Delta)\, s \ldots s\,\exp(t_k\Delta)\, dt\ , \label{duhamel}
\ee
where $\Delta_k= \{(t_0,\ldots,t_k)|\sum_{i=0}^k t_i=1\}$ is the standard $k$-simplex and $dt = dt_0\ldots dt_{k-1}$. 
\end{lemma}
{\it Proof:} Since the exponential of a generalized Laplacian is defined by its formal power series, the identity (\ref{duhamel}) which holds at a formal level makes sense in $\Sc(M)$. Indeed $s\in\Dc^0_k$ carries at least one power of $\eps$, so that the product $\exp(t_0\Delta) s \exp(t_1\Delta) s \ldots s\exp(t_k\Delta)$ is in $\Sc_k(M)$, and its expansion in powers of $\eps$ has polynomial coefficients with respect to $(t_0,\ldots,t_k)$. Hence the integral over the simplex $\Delta_k$ gives a well-defined element of $\Sc_k(M)$, and the infinite sum over $k$ converges in $\Sc(M)$. \cqfd\\

\noindent Notice that the Duhamel formula can be rewritten by means of the automorphism group $\si_{\Delta}$ as follows:
\be
\exp(\Delta+s) = \sum_{k=0}^{\infty} \int_{\Delta_k} \si^{t_0}_{\Delta}(s)  \si^{t_0+t_1}_{\Delta}(s) \ldots \si^{t_0+\ldots +t_{k-1}}_{\Delta}(s)\,\exp(\Delta) \, dt
\ee
Fix a generalized Laplacian $\Delta$ and consider the following vector subspace of $\Sc(M)$:
\be
\Tc(M) = \Dc(M) \exp\Delta
\ee
\begin{proposition}
$\Tc(M)$ is a $\Dc(M)$-bimodule and does not depend on the choice of generalized Laplacian. We call $\Tc(M)$ the bimodule of \emph{trace-class operators}. 
\end{proposition}
{\it Proof:}  $\Tc(M)$ is clearly a left $\Dc(M)$-module. Moreover by Lemma \ref{lsi}, one has $\Dc(M) \exp(\Delta) \Dc(M)= \Dc(M) \si^1_{\Delta}(\Dc(M))\exp\Delta= \Dc(M) \exp\Delta$ hence $\Tc(M)$ is a right $\Dc(M)$-module. Further on, if $\Delta$ and $\Delta'$ are two Laplacians, then $\Delta'=\Delta+s$ with $s\in \Dc^0_1(M)$. We know that $\si^t_{\Delta}(s)\in \Dc^0_1(M)$ for any $t\in\rr$, so the series 
$$
S=\sum_{k=0}^{\infty} \int_{\Delta_k} \si^{t_0}_{\Delta}(s)  \si^{t_0+t_1}_{\Delta}(s) \ldots \si^{t_0+\ldots +t_{k-1}}_{\Delta}(s)\, dt
$$
converges in $\Dc^0(M)$. Hence $\exp\Delta' = S \exp\Delta$ by the Duhamel formula. In the same way $\exp\Delta = S' \exp\Delta'$. Therefore $\Dc(M) \exp \Delta'=\Dc(M) \exp \Delta$, and $\Tc(M)$ does not depend on the choice of generalized Laplacian. \cqfd\\

\noindent $\Tc(M)$ is not a subalgebra of $\Sc(M)$. For example the product $\exp(\Delta) \exp(\Delta) = \exp(2\Delta)$ does not belong to the space of trace-class operators.

\section{Canonical trace}\label{sct}

Let $M$ be a closed manifold. The Wodzicki residue (\cite{W87}) is a canonical trace on the algebra of classical pseudodifferential operators $\CL(M)$. It is in fact the unique trace (up to a numerical factor) on $\CL(M)$ when the manifold has dimension $n>1$. The Wodzicki residue vanishes on $\CL^m(M)$ whenever $m<-n$, hence vanishes on the ideal of smoothing operators $\L^{-\infty}(M)$, so that it is really a trace on the algebra of formal symbols $\CS(M)$. Wodzicki gives a concrete formula for the residue of a symbol $a\in\CS^m(M)$ in terms of its expansion $a(x,p)=\sum_{j} a_{m-j}(x,p)$ in a local system of canonical coordinates over an open subset $U\subset M$. Let $\om=dp_i\wedge dx^i$ be the symplectic two-form on the cotangent bundle $T^*U\subset T^*M$ (summation over repeated indices). Then $T^*U$ is canonically oriented by the volume form $\om^n/n!= dp_1\wedge dx^1 \ldots dp_n\wedge dx^n$. The cosphere bundle $S^*U$ inherits this orientation. The Wodzicki residue of a symbol $a(x,p)$ with compact $x$-support over $U$ is the integral of a $(2n-1)$-form 
\be
\bint a = \frac{1}{(2\pi)^n} \int_{S^*U} \iota(L) \cdot \Big( a_{-n}(x,p)\  \frac{\om^n}{n!} \Big)\ , \label{wod}
\ee
where $a_{-n}$ is the degree $-n$ component of the symbol and $L=p_i\frac{\d}{\d p_i}$ is the fundamental vector field on $T^*U$ . We can write
$$
\iota(L) \cdot \frac{\om^n}{n!} = (\iota(L) \cdot \om) \wedge \frac{\om^{n-1}}{(n-1)!} = \frac{\eta\wedge\om^{n-1}}{(n-1)!} 
$$
where $\eta=p_idx^i$ is the canonical one-form on $T^*U$. It is non-trivial to check that the Wodzicki residue is a trace and does not depend on the choice of coordinate system. Hence such expressions can be patched together using a partition of unity, allowing to define the residue of a symbol $a$ with arbitrary support. If $E$ is a ($\zz_2$-graded) complex vector bundle over $M$, one defines analogously the Wodzicki residue as a (graded) trace on the algebra $\CS(M,E)$: at each point $(x,p)$ the symbol $a_{-n}(x,p)$ is now a endomorphism acting on the fibre $E_x$, so (\ref{wod}) has to be modified according to
\be
\bint a = \frac{1}{(2\pi)^n} \int_{S^*U} \iota(L) \cdot \Big( \tr_s\big(a_{-n}(x,p)\big)\ \frac{\om^n}{n!} \Big) \ ,
\ee
where $\tr_s$ is the (graded) trace of endomorphisms. We focus on $E=\Lambda T^*_{\cc}M$. In a local coordinate system over $U$ we know that a basis of sections of $\End(E)$ is provided by all products of $\psi^i$ or $\psib_j$, $i,j=1,\ldots,n$ among themselves, taking the Clifford relations (\ref{clif}) into account. A symbol $a\in \CS(U,E)$ may thus be decomposed into a finite sum over multi-indices $\eta = (\eta_1 , \ldots , \eta_n)$, $\theta = (\theta_1 , \ldots , \theta_n)$,
\be
a(x,p,\psi,\psib) = \sum_{\eta,\theta} a^{\eta,\theta}(x,p) \,\psi^{\eta} \psib^{\theta}\ ,\label{sum}
\ee
where the coefficients $a^{\eta,\theta}$ are functions of $x$ and $p$ only. It is easy to see that the graded trace of endomorphisms, which acts on polynomials $\psi^{\eta}\psib^{\theta}$, vanishes whenever $(|\eta|,|\theta|)\neq (n,n)$ and is normalized as follows on the polynomial of highest weight:
\be
\tr_s(\psi^1\ldots\psi^n \psib_{n}\ldots\psib_{1}) = (-1)^n \ .\label{super}
\ee
An equivalent normalization is $\tr_s(\Pi)=1$ where $\Pi=\psib_1\psi^1\ldots \psib_n\psi^n$ is the projection operator from the space of differential forms $\Om^*(U)$ to the subspace of scalar functions $\Om^0(U)$. \\

In section \ref{sbfs} we introduced the algebra $\Sc(M)$ acting on the space of formal power series $\CS(M,E)[[\eps]]$, its subalgebra $\Dc(M)\subset \Sc(M)$, and the $\Dc(M)$-bimodule of trace-class operators $\Tc(M)\subset \Sc(M)$. By means of the Wodzicki residue, our goal now is to construct a graded trace on $\Tc(M)$, that is, a linear map $\Tc(M)\to \cc$ vanishing on the subspace of graded commutators $[\Dc(M),\Tc(M)]$. We start by doing this locally on an open subset $U\subset M$. Choose a coordinate system $(x,p)$ over $U$ and fix the canonical ``flat" Laplacian $\Delta = \i\eps \frac{\d}{\d x^i}\frac{\d}{\d p_i}$. For all multi-indices $\al$ and $\beta$ set
\be
\langle \d_x^{\al}\d_p^{\beta} \exp\Delta \rangle = \left. \d_x^{\al}\d_p^{\beta} \cdot \exp\Big(\frac{\i}{\eps} (p_i-q_i)(x^i-y^i)\Big) \right|_{\substack{x=y \\ p=q}} \label{cont}
\ee
For example one has
$$
\langle \exp\Delta \rangle = 1\ ,\quad \langle \d_{x^i} \exp\Delta\rangle = 0 = \langle \d_{p_j} \exp\Delta\rangle\ ,\quad \langle \d_{x^i}\d_{p_j}\exp\Delta \rangle = \frac{\i}{\eps}\, \delta^j_i
$$
and more generally with a polynomial $\d_x^{\al}\d_p^{\beta}$ the formula involves all possible contractions between $\d_x$ and $\d_p$. In particular
$$
\langle \d_{x^i}\d_{x^j}\d_{p_k}\d_{p_l} \exp\Delta \rangle = \Big(\frac{\i}{\eps}\Big)^2 (\delta^k_i\delta^l_j + \delta^l_i\delta^k_j) \ . 
$$
Notice that $\langle \d_x^{\al}\d_p^{\beta} \exp\Delta\rangle$ vanishes unless $|\al|=|\beta|$. We define similarly a contraction map for the polynomials in the odd variables $\psi_R,\psib_R$. If $(\psi^{\eta}\psib^{\theta})_R$ is a generic product with multi-indices $\eta,\theta$ set
\be
\langle (\psi^{\eta}\psib^{\theta})_R \rangle = (-1)^n \tr_s(\psi^{\eta} \psib^{\theta})\ .
\ee
Hence from the normalization (\ref{super}) holds $\langle \psi^1\ldots\psi^n \psib_{n}\ldots\psib_{1} \rangle=1$, and the contraction vanishes on polynomials of lower degree. The even and odd contractions assemble in a linear map
\be
\langle\!\langle \ \rangle\!\rangle : \Tc(U) \to \CS(U,E)[[\eps]]
\ee
defined as follows. Let $s=\sum_{k=0}^{\infty}s_k\eps^k$ belong to $\Dc^m(U)$, so that $s\exp\Delta$ is a generic element of $\Tc(U)$. We can write, for all components $s_k\in\Lc(U)$,
$$
s_k = \sum_{|\al|=0}^{k} \sum_{|\beta|=0}^{\infty} \sum_{|\eta|=0}^n \sum_{|\theta|=0}^n (s_{k,\al,\beta,\eta,\theta})_L \, (\psi^{\eta}\psib^{\theta})_R \, \d_x^{\al} \, \d_p^{\beta} 
$$
with $s_{k,\al,\beta,\eta,\theta}\in \CS(U,E)$ a symbol of order $\leq m+(k+|\beta|-3|\al|)/2$. Set  
$$
\langle\!\langle s_k \exp\Delta\rangle\!\rangle = \sum_{|\al|=0}^{k} \sum_{|\beta|=0}^{\infty} \sum_{|\eta|=0}^n \sum_{|\theta|=0}^n   s_{k,\al,\beta,\eta,\theta} \, \langle(\psi^{\eta} \psib^{\theta})_R \rangle  \langle \d_x^{\al} \, \d_p^{\beta} \exp\Delta \rangle \ .
$$
Observe that the sum over $\al$ is finite, as is the sum over $\beta$ because of the contractions $\langle \d_x^{\al} \, \d_p^{\beta} \exp\Delta\rangle$. Hence $\langle\!\langle s_k \exp\Delta\rangle\!\rangle$ is a polynomial of degree at most $k$ in the indeterminate $\eps^{-1}$, with coefficients in $\CS(U,E)$. Consequently $\langle\!\langle s_k \exp\Delta\rangle\!\rangle \eps^k$ is a polynomial in $\eps$ of degree at most $k$, with coefficients in $\CS(U,E)$. However it is not at all obvious that the sum
\be
\langle\!\langle s \exp\Delta\rangle\!\rangle = \sum_{k=0}^{\infty}\langle\!\langle s_k \exp\Delta\rangle\!\rangle \eps^k
\ee
makes sense even in the space of formal series $\CS(U,E)[[\eps]]$. The completeness of the space of symbols is an essential ingredient of the following lemma. 

\begin{lemma}
$\langle\!\langle s \exp\Delta\rangle\!\rangle$ is a well-defined element of $\CS(U,E)[[\eps]]$ for any $s\in \Dc(U)$.  
\end{lemma}
{\it Proof:} Let $s\in \Dc^m(U)$. For each power $l\in\nn$, we have to show that the coefficient of $\eps^l$ in the formal series 
$$
\langle\!\langle s \exp\Delta\rangle\!\rangle = \sum_{k=0}^{\infty}\sum_{|\al|=0}^{k} \sum_{|\beta|=0}^{\infty} \sum_{|\eta|=0}^n \sum_{|\theta|=0}^n   s_{k,\al,\beta,\eta,\theta} \, \langle(\psi^{\eta} \psib^{\theta})_R\rangle \langle \d_x^{\al} \, \d_p^{\beta} \exp\Delta \rangle \eps^k
$$
is a well-defined element of $\CS(U,E)$. The contraction $\langle \d_x^{\al} \, \d_p^{\beta} \exp\Delta \rangle$ forces $|\beta|=|\al|$, hence the symbol $s_{k,\al,\beta,\eta,\theta}$ has order $\leq m+k/2 -|\al|$. Moreover $\langle \d_x^{\al} \, \d_p^{\beta} \exp\Delta \rangle$ brings a factor $\eps^{-|\al|}$. It follows that for fixed $l\in\nn$, the coefficient of $\eps^l$ in the above series is proportional to 
$$
a_l=\sum_{k=0}^{\infty}\sum_{|\al|=k-l}  a_{k,\al} 
$$
where $a_{k,\al}$ is a symbol of order $\leq m+k/2 -|\al|=m+l-k/2$. Since $m$ and $l$ are fixed, the order of $a_{k,\al}$ is a strictly decreasing function of $k$, hence $a_l$ converges in $\CS(U,E)$. \cqfd\\

\noindent Let $\Dc_c(U)\subset \Dc(U)$ and $\Tc_c(U)\subset \Tc(U)$ be the subspaces of operators with compact $x$-support on $U$. Any element of $\Tc_c(U)$ reads $s\exp\Delta$ for some $s\in \Dc_c(U)$, and $\Tc_c(U)$ is a $\Dc(U)$-bimodule.

\begin{lemma}
Let $\langle\!\langle s \exp\Delta\rangle\!\rangle [n] \in \CS(U,E)$ be the coefficient of $\eps^n$, $n=\dim M$, in the formal series $\langle\!\langle s \exp\Delta\rangle\!\rangle$. The linear map $\Tr_s^U:\Tc_c(U)\to \cc$ defined by
\be
\Tr_s^U(s\exp\Delta) = \bint \langle\!\langle s \exp\Delta\rangle\!\rangle [n]\ , \qquad \forall \ s\in \Dc_c(U)\ , 
\ee
is a graded trace on the space of compactly-supported trace-class operators viewed as a $\Dc(U)$-bimodule. 
\end{lemma}
{\it Proof:} In fact we will show that the map $\Tc_c(U)\to \cc[[\eps]]$ defined by 
$$
s\exp\Delta \mapsto \bint \langle\!\langle s \exp\Delta\rangle\!\rangle
$$
is a graded trace. Selecting the coefficient of $\eps^n$ then yields $\Tr_s^U$. By linearity it is sufficient to check the trace property on operators $s\in \Dc(U)$ which depend \emph{polynomially} on $\eps$ and the partial derivatives $\d_x$ and $\d_p$. So let $s=s_k\eps^k$,
$$
s_k= a_L (\psi^{\eta}\psib^{\theta})_R \d_x^{\al}\d_p^{\beta}
$$
be such an operator, for some $a\in \CS(U,E)$ and multi-indices $\al,\beta,\eta,\theta$. It is enough to show that 
\be
\bint \langle\!\langle (s \exp\Delta) s'\rangle\!\rangle = \pm \bint \langle\!\langle s' s \exp\Delta \rangle\!\rangle \label{trace}
\ee
in the following cases: $s'=b_L$ for a symbol $b\in \CS(U,E)$, or $s'=\d_x$, $\d_p$, $\psi_R$, $\psib_R$. The sign must be $-$ if $s$ and $s'$ are both odd, $+$ otherwise. Since the contraction map involves the supertrace on the Clifford algebra generated by $\psi_R,\psib_R$, (\ref{trace}) is obvious when $s'=\psi_R$ or $\psib_R$. Then for $s'=\frac{\d}{\d x^i}$ one has
$$
\bint \langle\!\langle [s', s \exp\Delta] \rangle\!\rangle = \bint \frac{\d a}{\d x^i} \langle(\psi^{\eta}\psib^{\theta})_R\rangle \langle \d_x^{\al}\d_p^{\beta}\exp\Delta \rangle  \eps^k
$$
The Wodzicki residue vanishes on the derivative $\d a/\d x^i$, hence (\ref{trace}) is verified. The case $s'=\frac{\d}{\d p_i}$ is similar. It remains to deal with the case $s'=b_L$ for a symbol $b$. If $F(\d_x,\d_p)$ is any formal power series with respect to the variables $X=\d_x$ and $P=\d_p$, one has the identity
$$
F(\d_x,\d_p)\circ b_L = \sum_{|\gamma|=0}^{\infty} \sum_{|\delta|=0}^{\infty} \frac{1}{\gamma! \delta!} \,  (\d_x^{\gamma} \d_p^{\delta} b)_L \, \d_X^{\gamma} \d_P^{\delta}F(\d_x,\d_p) \ .
$$
Applying this to the series $F(\d_x,\d_p) = \d_x^{\al}\d_p^{\beta} \exp\Delta$ one gets
$$
\langle (a_L \d_x^{\al}\d_p^{\beta} \exp\Delta) b_L \rangle = \sum_{|\gamma|=0}^{\infty} \sum_{|\delta|=0}^{\infty} \frac{1}{\gamma! \delta!} \, \langle (a \d_x^{\gamma} \d_p^{\delta} b)_L \,  \d_X^{\gamma} \d_P^{\delta} (\d_x^{\al}\d_p^{\beta} \exp\Delta)\rangle\ .
$$
But the contraction map vanishes on a derivative $\d_XF(\d_x,\d_p)$ or $\d_PF(\d_x,\d_p)$. This only selects the terms $|\gamma|=|\delta|=0$:
$$
\langle (a_L \d_x^{\al}\d_p^{\beta} \exp\Delta) b_L \rangle = \langle (ab)_L \, \d_x^{\al}\d_p^{\beta} \exp\Delta \rangle = ab \,\langle  \d_x^{\al}\d_p^{\beta} \exp\Delta \rangle \ .
$$
Finally, the Wodzicki residue is a trace on the algebra of (compactly supported) symbols, hence
$$
\bint \langle (a_L \d_x^{\al}\d_p^{\beta} \exp\Delta) b_L \rangle = \bint ba \,\langle  \d_x^{\al}\d_p^{\beta} \exp\Delta \rangle = \bint \langle b_L a_L \d_x^{\al}\d_p^{\beta} \exp\Delta \rangle\ .
$$
This shows that (\ref{trace}) is verified for $s'=b_L$ as well. \cqfd\\

\begin{proposition}
The map $\Tr_s^U$ does not depend on the choice of coordinate system $(x,p)$ over $T^*U$. Hence using a partition of unity relative to an open covering of $M$, such maps can be patched together, giving rise to a \emph{canonical} graded trace
\be
\Tr_s: \Tc(M)\to \cc
\ee
on the $\Dc(M)$-bimodule of trace-class operators.
\end{proposition}
{\it Proof:} First observe that under an \emph{affine} change of coordinates $x^i\mapsto y^i$, $p_i\mapsto \frac{\d x^j}{\d y^i} p_j$, the flat laplacian $\Delta = \i\eps \frac{\d}{\d x^i}\frac{\d}{\d p_i}$ is invariant, as well as Eq. (\ref{cont}). It follows that the contraction map $\Tc(U)\to  \CS(U,E)[[\eps]]$ is equivariant under affine transformations. Since the Wodzicki residue is also invariant, it follows that the trace $\Tr^U_s$ is invariant under affine transformations.\\
Now let $\gamma$ be any (smooth) change of coordinates. By linearity it is enough to show that, if $t\in\Tc(U)$ has support in an arbitrary small neighborhood of a point $x_0\in U$, then $\Tr^U_s(t)=\Tr^U_s(\gamma(t))$. After composition with an appropriate affine transformation, we can even suppose that $\gamma$ leaves the point $x_0$ and its tangent space $T_{x_0}U$ fixed. Then there exists a small neighborhood $V$ of $x_0$ such that the restriction of $\gamma$ to the domain $V$ is a diffeomorphism homotopic to identity. Hence we only need to show that the trace $\Tr_s^U$ is invariant under infinitesimal transformations induced by vector fields on $U$. This follows from the fact that such transformations are given by commutators. Indeed let $X=X^i\frac{\d}{\d x^i} \in \Vect(U)$ be any smooth vector field and consider the following symbol $L_X\in\PS(U,E)$:
$$
L_X=\i X^jp_j + \frac{\d X^j}{\d x^k} \psi^k\psib_j \ .
$$
As an operator on the smooth sections of $E$ (that is, the differential forms over $U$), $L_X$ corresponds to the Lie derivative along $X$. One easily checks that its induced action on the generators of the algebra $\CS(U,E)$ reads
\beq
{[L_X, x^i]} &=& X^i\ , \qquad  [L_X,p_i]= -\frac{\d X^j}{\d x^i}p_j +\i \frac{\d^2 X^j}{\d x^i\d x^k}\psi^k\psib_j\ , \non\\
{[L_X,\psi^i]} &=& \frac{\d X^i}{\d x^k}\psi^k\ ,\qquad [L_X,\psib_i] = -\frac{\d X^j}{\d x^i}\psib_j\ ,\non
\eeq
which are the correct transformation laws. Further on, the induced action on the algebra $\Sc(U)$ is given by the commutator with $(L_X)_L+(L_X)_R \in \Lc(U)$. Restricting this action to the subspace of trace-class operators $\Tc(U)$ shows that the trace $\Tr^U_s$ vanishes on Lie derivatives. \cqfd\\

We end this section with a useful formula in local coordinates $(x,p)$ over $U\subset M$. Let $R=(R^i_j)$ be an $n\times n$ matrix with entries in $\cc[[\eps]]$. We suppose that $R$ has no term of degree zero with respect to $\eps$. Hence the Todd series of $R$ and its determinant are well-defined as formal power series in $M_n(\cc[[\eps]])$ and $\cc[[\eps]]$ respectively:
\be
\frac{R}{e^R-1} = 1 - \frac{1}{2}R + \frac{1}{12}R^2 + \ldots \ ,\qquad \Td(R) = \det\left( \frac{R}{e^R-1} \right)\ .
\ee
We consider the operator $s= p_{iL}R^i_j\d_{p_j}=p_L\cdot R \cdot \d_p$ as a formal perturbation of the flat Laplacian $\Delta = \i\eps \d_x\cdot \d_p$. Note however that $\Delta+s$ is not a generalized Laplacian. Then by the Duhamel formula,
$$
\exp(\Delta+p_L\cdot R \cdot \d_p) = \sum_{k=0}^{\infty} \int_{\Delta_k} \big( \si^{t_0}_{\Delta}(s)  \si^{t_0+t_1}_{\Delta}(s) \ldots \si^{t_0+\ldots +t_{k-1}}_{\Delta}(s)\, \exp\Delta\big) dt
$$
is a well-defined element of $\Tc(U)$. There is an explicit formula computing the contraction of this series with an arbitrary polynomial in the derivatives $\d_x$ and $\d_p$:

\begin{lemma}\label{luse}
For any multi-indices $\al$ and $\beta$ holds
\be
\langle \d_x^{\al}\d_p^{\beta} \exp(\Delta+ p_L\cdot R \cdot \d_p) \rangle = \Td(R) \, s(R,p) 
\ee
where the symbol $s(R,p)\in \CS(U,E)[[\eps]]$ is a polynomial in $p$:
$$
s(R,p)= \left. \d_x^{\al}\d_p^{\beta} \exp\Big(\frac{\i}{\eps} q\cdot R \cdot (x-y) + \frac{\i}{\eps} (p-q)\cdot \frac{R}{1-e^{-R}} \cdot (x-y)\Big) \right|_{\substack{x=y \\ p=q}} 
$$
\end{lemma}
{\it Proof:} The operator $\exp(\Delta+p_L\cdot R\cdot \d_p)\exp(-\Delta)$ can be expanded as a formal power series in $R$, whose coefficients depend polynomially on $p_L$ and the partial derivatives $\d_x,\d_p$. Thus one has
$$
\langle \d_x^{\al}\d_p^{\beta} \exp(\Delta+ p_L\cdot R\cdot\d_p) \rangle =  \d_x^{\al}\d_p^{\beta} H_{\eps}(R,x,y,p,q) \Big|_{\substack{x=y \\ p=q}} 
$$
where $H_{\eps}(R,x,y,p,q)=\exp(\Delta+p\cdot R\cdot \d_p)\exp(-\Delta) \big( \exp\big( \frac{\i}{\eps} (p-q)\cdot(x-y)\big)\big)$. We introduce a deformation parameter $t\in[0,1]$ and replace $R$ by $tR$. The function $H_{\eps}(tR,x,y,p,q)$ is viewed as a formal power series in $t$. For $t=0$ it reduces to
$$
H_{\eps}(0,x,y,p,q) = \exp\Big(\frac{\i}{\eps} (p-q)\cdot(x-y)\Big)\ .
$$
We are going to show that $H_{\eps}(tR,x,y,p,q)$ fulfills a differential equation of first order with respect to $t$. One has
$$
\frac{\d}{\d t} \exp(\Delta+p\cdot tR\cdot \d_p) = \sum_{n=0}^{\infty}\frac{1}{(n+1)!} (p\cdot R\cdot\d_p)^{(n)} \exp(\Delta+p\cdot tR\cdot \d_p)
$$
where the superscript $^{(n)}$ denotes the derivation $X\mapsto [\Delta+p\cdot tR\cdot \d_p, X]$ applied $n$ times. Hence $(p\cdot R\cdot\d_p)^{(1)}= [\Delta,p\cdot R\cdot\d_p]=\i \eps \d_x\cdot R\cdot\d_p = \i\eps R^i_j\frac{\d}{\d x^i}\frac{\d}{\d p_j}$. Furthermore 
$$
(p\cdot R\cdot\d_p)^{(n)} = [p\cdot tR\cdot\d_p , (p\cdot R\cdot\d_p)^{(n-1)}] = \i\eps (-t)^{n-1} \d_x\cdot R^n\cdot\d_p
$$
for all $n\geq 2$. Hence we can write
\beq
\lefteqn{\frac{\d}{\d t} \exp(\Delta+p\cdot tR\cdot \d_p)}\non\\ 
&&\qquad = \Big(p\cdot R\cdot\d_p - \i\eps \sum_{n=1}^{\infty} \d_x \cdot \frac{(-tR)^n }{t(n+1)!} \cdot\d_p \Big) \exp(\Delta+p\cdot tR\cdot \d_p) \non\\
&&\qquad = \Big( p\cdot R\cdot\d_p + t^{-1} \Delta + \i\eps\, \d_x \cdot \frac{e^{-tR}-1 }{t^2R} \cdot\d_p \Big) \exp(\Delta+p\cdot tR\cdot \d_p) \non
\eeq
It is important to note that, by construction, the $t$-expansion of the differential operator $p\cdot R\cdot\d_p + t^{-1} \Delta + \i\eps\, \d_x \cdot \frac{e^{-tR}-1 }{t^2R} \cdot\d_p$ only involves non-negative powers of $t$. Hence the function $H_{\eps}$ is a solution of the differential equation 
$$
\Big( -\frac{\d}{\d t} + p\cdot R\cdot\d_p + t^{-1} \Delta + \i\eps\, \d_x \cdot \frac{e^{-tR}-1 }{t^2R} \cdot\d_p \Big) \, H_{\eps}(tR,x,y,p,q) = 0 
$$
and is uniquely specified, as a formal power series in $t$, by its value at $t=0$. A routine computation shows that the Ansatz
$$
H_{\eps}(tR,x,y,p,q) = \Td(tR) \exp\Big(\frac{\i}{\eps} q\cdot tR \cdot (x-y) + \frac{\i}{\eps} (p-q)\cdot \frac{tR}{1- e^{-tR}} \cdot (x-y)\Big)
$$
is this unique solution. \cqfd \\

\noindent Let us apply this lemma in some particular cases. One has
\beq
\langle \exp(\Delta+ p_L\cdot R \cdot \d_p) \rangle
&=& \Td(R) \\
\Big\langle \frac{\d}{\d x^i} \exp(\Delta+p_L\cdot R \cdot \d_p) \Big\rangle
&=& \frac{\i}{\eps}\, \Td(R) \, (p\cdot R)_i \non
\eeq 
Observe that the right-hand-side of the second equation contains no negative power of $\eps$ because $R$ brings at least one factor $\eps$. One thus gets the identity
\be
\Big\langle \Big(\i\eps \frac{\d}{\d x^i} + (p_L\cdot R)_i \Big) \exp(\Delta+ p_L\cdot R \cdot \d_p) \Big\rangle = 0\ .
\ee
More generally for any multi-index $\al=(\al_1,\ldots,\al_n)$:
\be
\big\langle \big(\i\eps\d_x+p_L\cdot R \big)^{\al} \exp(\Delta+ p_L\cdot R \cdot \d_p) \big\rangle = 0\ . \label{iden}
\ee

\section{Dirac operators}\label{sdo}

Let $M$ be an $n$-dimensional manifold and $E=\Lambda T^*_{\cc}M$. The space of smooth sections of $E$ is isomorphic to the space $\Om^*(M)$ of complex differential forms over $M$. The exterior multiplication of $\Om^*(M)$ on the sections of $E$ (from the left) gives rise to an homomorphism of algebras
\be
\mu: \Om^*(M)\to \PS^0(M,E)\ .
\ee
Remark that the algebra $\PS^0(M,E)$ of differential operators of order zero is isomorphic to the algebra of smooth sections of the endomorphism bundle $\End(E)$. The map $\mu$ is injective. In a local coordinate system $(x^1,\ldots,x^n)$ over $U\subset M$, the image of a $k$-form $\al= \al_{i_1\ldots i_k}(x)dx^{i_1}\wedge\ldots\wedge dx^{i_k}$ is the endomorphism $\mu(\al)= \al_{i_1\ldots i_k}(x)\psi^{i_1}\ldots \psi^{i_k}$. Also, the operation of interior multiplication by vector fields on the sections of $E$ gives rise to an injective linear map
\be
\iota: \Vect(M)\to \PS^0(M,E)\ .
\ee
In local coordinates the image of a vector field $X= X^i(x)\frac{\d}{\d x^i}$ is the endomorphism $\iota(X)= X^i(x)\psib_i$. In the sequel we consider $\Om^*(M)$ and $\Vect(M)$ as subspaces of the algebra of differential operators $\PS(M,E)$. Finally, we introduce another subspace $\SPS^1(M,E)\subset\PS^1(M,E)$ of differential operators, characterized by their expression in any local coordinate system over $U$ as follows:
\be
a\in \SPS^1(U,E) \Leftrightarrow a(x,p)= a^i(x)p_i + a^i_j(x) \psi^j\psib_i + b(x)
\ee
where $a_i,a^i_j,b \in \Om^0(U)$ are scalar functions. $\SPS^1(M)$ is the space of differential operators of order one, even parity, and \emph{scalar} leading symbol. This definition is coordinate-independent, because under a coordinate change $x^i\mapsto y^i$ one has $p_i\mapsto \frac{\d x^j}{\d y^i}p_j -\i \frac{\d}{\d x^k}\frac{\d x^j}{\d y^i}\psi^k\psib_j$, $\psi^i \mapsto \frac{\d y^i}{\d x^j}\psi^j$ and $\psib_i \mapsto \frac{\d x^j}{\d y^i}\psib_j$. Moreover, one easily checks that the commutator $[\SPS^1(M,E),\SPS^1(M,E)]$ is again in $\SPS^1(M,E)$. Thus $\SPS^1(M,E)$ is a Lie algebra. \\
As before we denote $\Lc(M)$ the subalgebra of linear operators on the vector space $\CS(M,E)$, generated by left multiplication by $\CS(M,E)$ and right multiplication by $\PS(M,E)$. In other words $\Lc(M)=\CS(M,E)_L\PS(M,E)_R$. From the discussion above we can also form various subspaces of $\Lc(M)$, for instance $\SPS^1(M,E)_L\Om^1(M)_R$ or $\Om^0(M)_L\Vect(M)_R$. The latter operators are easy to characterize in local coordinates:
\beq
s\in \SPS^1(U,E)_L\Om^1(U)_R &\Leftrightarrow& s = \sum_{|\al|=0}^{\infty} (s^k_{\al i}p_k + s^k_{\al ij}\psi^j\psib_k + s_{\al i})_L\psi^i_R \d_p^{\al}   \non\\
r\in \Om^0(U)_L\Vect(U)_R  &\Leftrightarrow& r = \sum_{|\al|=0}^{\infty} (r^i_{\al})_L \psib_{iR} \d_p^{\al}
\eeq
for some scalar functions $s^k_{\al i}, s^k_{\al ij}, s_{\al i}, r^i_{\al} \in \Om^0(U)$. From these expressions it is clear that $\SPS^1(M,E)_L\Om^1(M)_R \subset \Dc^1_0(M)$ and $\Om^0(M)_L\Vect(M)_R \subset \Dc^0_0(M)$. We are now ready to define Dirac operators as particular elements of $\Dc(M)$.

\begin{definition}\label{ddirac}
Suppose that $\nabla\in \Lc(M)$ and $\nablab\in\Lc(M)$ are odd operators on $\CS(M,E)$ such that, in any local coordinate system over $U\subset M$, 
\be
\nabla = \psi^i_R\frac{\d}{\d x^i} + s \ ,\qquad \nablab = \psib_{iR}\frac{\d}{\d p_i} + r\ ,
\ee
with $s\in \SPS^1(U,E)_L\Om^1(U)_R$ and $r\in \Om^0(U)_L\Vect(U)_R \cap \Dc^{-1}_0(U)$. The sum 
\be
D=\i\eps\nabla + \nablab\ \in \ \Dc^1_1(M) + \Dc^{-1/2}_0(M)
\ee
is called a \emph{generalized Dirac operator} on $M$.
\end{definition}
Hence the operator $\i\eps\nabla$ is locally the sum of its leading part $\i\eps\psi^i_R\frac{\d}{\d x^i}\in \Dc^{1}_1(U)$ and a perturbation term $\i\eps s\in \Dc^{1/2}_1(U)$. Similarly $\nablab$ is locally the sum of its leading part $\psib_{iR}\frac{\d}{\d p_i} \in \Dc^{-1/2}_0(U)$ and a perturbation $$
r= \sum_{|\al|=2}^{\infty} (r^i_{\al})_L \psib_{iR} \d_p^{\al}\ \in \Dc^{-1}_0(U)\ .
$$
In fact $\psib_{iR}\frac{\d}{\d p_i}=\i \psib_{iR}(x^i_R-x^i_L)$ so that $\nablab \in \Om^0(M)_L\Vect(M)_R \cap \Dc^{-1/2}_0(M)$. \\

The existence of $\nabla$ and $\nablab$ as global operators on $M$ is not a priori obvious. In order to understand the above definition, we first examine the behaviour of $\psi^i_R\frac{\d}{\d x^i}$ under a coordinate change $x^i\mapsto\gamma(x^i)=y^i$. One has
$$
\gamma\Big( \psi^i_R\frac{\d}{\d x^i} \Big) = \gamma(\psi^i_R)\gamma(\i p_{iL}-\i p_{iR}) = \i \gamma(\psi^i)_R\gamma(p_i)_L -\i \gamma(p_i\psi^i)_R\ .
$$
$\i p_i\psi^i$ is the symbol of the exterior derivative $d$ on differential forms (see Example \ref{ecan}), hence it is ivariant under coordinate change. This can also be checked by direct computation, with $\gamma(p_i)=\frac{\d x^j}{\d y^i}p_j -\i \frac{\d}{\d x^k}\frac{\d x^j}{\d y^i}\psi^k\psib_j$ and $\gamma(\psi^i)=\frac{\d y^i}{\d x^l}\psi^l$. One thus has
$$
\gamma\Big( \psi^i_R\frac{\d}{\d x^i} \Big) = \i \Big(\frac{\d y^i}{\d x^l}\psi^l\Big)_{\!\! R} \Big(\frac{\d x^j}{\d y^i}p_j -\i \frac{\d}{\d x^k}\frac{\d x^j}{\d y^i}\psi^k\psib_j\Big)_{\!\! L} -\i (p_i\psi^i)_R\ .
$$
Write $-\i (p_i\psi^i)_R=\psi^i_R\frac{\d}{\d x^i} -\i p_{iL}\psi^i_R$ and use the commutation of left and right actions:
\be
\gamma\Big( \psi^i_R\frac{\d}{\d x^i} \Big) = \psi^i_R\frac{\d}{\d x^i} -\i p_{iL}\psi^i_R + \Big(\i \frac{\d x^j}{\d y^i}p_j + \frac{\d}{\d x^k}\frac{\d x^j}{\d y^i}\psi^k\psib_j\Big)_{\!\! L} \Big(\frac{\d y^i}{\d x^l}\psi^l\Big)_{\!\! R}\ . \label{change}
\ee
The right-hand side reads $\psi^i_R\frac{\d}{\d x^i} + s$ with $s\in \SPS^1(U,E)_L\Om^1(U)_R$. We can choose a partition of unity $(c_I)$ relative to an atlas $(U_I,x_I)$ of $M$ and define a global operator $\nabla$ by:
\be
\nabla = \sum_{I} (c_I)_L (\psi^i_I)_R \frac{\d}{\d x^i_I} \ .
\ee
Then (\ref{change}) shows that $\nabla$ has the required form in any local coordinate system. In order to build a global operator $\nablab$, we proceed analogously and examine how the local operator $\psib_{iR}\frac{\d}{\d p_i}$ transforms under coordinate change. One has
$$
\gamma\Big(\psib_{iR}\frac{\d}{\d p_i}\Big)=\gamma(\psib_{iR}) \gamma(\i(x^i_R-x^i_L)) = \i \Big(\frac{\d x^j}{\d y^i}\psib_j\Big)_{\!\!R} (y^i_R-y^i_L)\ .
$$
This still belongs to the subspace $\Om^0(U)_L\Vect(U)_R \cap \Dc^{-1/2}_0(U)$. Then use the expansion $y^i_R=\sum_{|\al|=0}^{\infty} \frac{(-\i)^{|\al|}}{\al !} \, (\d_x^{\al} y^i)_L \d_p^{\al}$:
$$
\gamma\Big(\psib_{iR}\frac{\d}{\d p_i}\Big) = \i \Big(\frac{\d x^j}{\d y^i}\psib_j\Big)_{\!\!R} \sum_{|\al|=1}^{\infty} \frac{(-\i)^{|\al|}}{\al !} \, (\d_x^{\al} y^i)_L \d_p^{\al} \ .
$$
For $|\al|\geq 2$ the terms of the series belong to $\Dc^{-1}_0(U)$. Keeping only the first term ($|\al|=1$) one gets
$$
\gamma\Big(\psib_{iR}\frac{\d}{\d p_i}\Big) \equiv \Big(\frac{\d x^j}{\d y^i}\psib_j\Big)_{\!\!R} \Big(\frac{\d y^i}{\d x^k}\Big)_{\!\!L} \frac{\d}{\d p_k} \mod \Dc^{-1}_0(U)\ .
$$
But $\big(\frac{\d x^j}{\d y^i}\big)_{\!R}= \sum_{|\al|=0}^{\infty} \frac{(-\i)^{|\al|}}{\al !} \, \big(\d_x^{\al} \frac{\d x^j}{\d y^i}\big)_{\! L} \d_p^{\al}$ equals $\big(\frac{\d x^j}{\d y^i}\big)_{\!L}$ modulo $\Dc^{-1/2}_0(U)$, so 
$$
\gamma\Big(\psib_{iR}\frac{\d}{\d p_i}\Big) \equiv \psib_{jR} \Big(\frac{\d x^j}{\d y^i}\frac{\d y^i}{\d x^k}\Big)_{\!\!L} \frac{\d}{\d p_k}\equiv \psib_{jR}\frac{\d}{\d p_j} \mod \Dc^{-1}_0(U)  \ .
$$
Hence we have $\gamma(\psib_{iR}\frac{\d}{\d p_i})=\psib_{iR}\frac{\d}{\d p_i}+r$ with $r\in \Om^0(U)_L\Vect(U)_R \cap \Dc^{-1}_0(U)$. As before we obtain a global operator $\nablab$ on $M$ by gluing local pieces $\psib_{iR}\frac{\d}{\d p_i}$ together by means of a partition of unity. This shows that generalized Dirac operators always exist.

\begin{proposition}\label{plap}
Let $D$ be a generalized Dirac operator on $M$. Then $-D^2$ is a generalized Laplacian (Definition \ref{dlap}). 
\end{proposition}
{\it Proof:} $D=\i\eps\nabla+\nablab$ so $D^2=\nablab^2 +\i\eps [\nabla,\nablab] -\eps^2\nabla^2$. Since $\nablab\in \Om^0(M)_L\Vect(M)_R$, one has $\nablab^2=0$ (indeed $\Om^0(M)$ is a commutative algebra, and vector fields anticommute in $\PS^0(M,E)$). Then choose a local coordinate system and write $\nabla=\psi^i_R\frac{\d}{\d x^i} +s$ with $s\in \SPS^1_L\Om^1_R$. One has $(\psi^i_R\frac{\d}{\d x^i})^2=0$ hence
$$
\nabla^2 = \Big(\psi^i_R\frac{\d}{\d x^i} +s\Big)^2 = [\psi^i_R\frac{\d}{\d x^i},s] + s^2\ .
$$
Recall that $\psi^i_R\frac{\d}{\d x^i}$ and $s$ are odd, so the commutator is taken in the graded sense. Let us decompose $s$ as a sum of generic elements $a_Lb_R$, with $a\in \SPS^1$ (even) and $b\in\Om^1$ (odd). Then $\psi^i$ and $b$ commute in the graded sense so
$$
[\psi^i_R\frac{\d}{\d x^i},a_Lb_R]= \psi^i_R[\frac{\d}{\d x^i},a_Lb_R] =  \Big(\frac{\d a}{\d x^i}\Big)_{\!\! L}\psi^i_R b_R +  a_L\psi^i_R\Big(\frac{\d b}{\d x^i}\Big)_{\!\! R} \ .
$$
Hence $[\psi^i_R\frac{\d}{\d x^i},s]\in \SPS^1_L\Om^2_R$. Then writing $s=\sum_{I}a^I_Lb^I_R$, one has
$$
s^2 = \sum_{I,J} a^I_La^J_Lb^I_R b^J_R = - \sum_{I,J} (a^Ia^J)_L(b^Jb^I)_R = - \frac{1}{2} \sum_{I,J} [a^I,a^J]_L(b^Jb^I)_R
$$
because $b^I$ and $b^J$ are anticommuting one-forms. $\SPS^1$ is a Lie algebra, hence $[a^I,a^J] \in \SPS^1$ and $s^2\in \SPS^1_L\Om^2_R$. This shows that $\nabla^2 \in \SPS^1_L\Om^2_R \subset \Dc^1_0$ and $-\eps^2\nabla^2 \in \Dc^0_2$. \\
Finally we compute the graded commutator $[\nabla,\nablab]$. Write $\nablab=\psib_{iR}\frac{\d}{\d p_i}+r$ with $r\in \Om^0_L\Vect_R \cap \Dc^{-1}_0$. One has $[\nabla, \psib_{iR}\frac{\d}{\d p_i}] = [\psi^j_R\frac{\d}{\d x^j}, \psib_{iR}\frac{\d}{\d p_i}] + [s, \psib_{iR}\frac{\d}{\d p_i}]$, where 
$$
[\psi^j_R\frac{\d}{\d x^j}, \psib_{iR}\frac{\d}{\d p_i}]=[\psi^j_R, \psib_{iR}]\frac{\d}{\d x^j}\frac{\d}{\d p_i}= - [\psib_i , \psi^j]_R\frac{\d}{\d x^j}\frac{\d}{\d p_i} = - \frac{\d}{\d x^i}\frac{\d}{\d p_i}
$$
As before decompose $s$ as a sum of generic elements $a_Lb_R\in \SPS^1_L\Om^1_R$.  Since $b=\sum_i b_i(x)\psi^i \in \Om^1$ does not depend on $p$, 
$$
[a_Lb_R, \psib_{iR}\frac{\d}{\d p_i}] = a_L[b_R, \psib_{iR}]\frac{\d}{\d p_i} + \psib_{iR} \Big( \frac{\d a}{\d p_i} \Big)_{\!\! L} b_R = - a_Lb_{iR}\frac{\d}{\d p_i} - \Big( \frac{\d a}{\d p_i} \Big)_{\!\! L} (b\psib_i)_R
$$
One has $a_Lb_{iR}\in \SPS^1_L\Om^0_R $, hence $a_Lb_{iR}\frac{\d}{\d p_i} \in \SPS^1_L\Om^0_R \cap \Dc^{1/2}_0$. Moreover $\big( \frac{\d a}{\d p_i} \big)_{\! L} (b\psib_i)_R \in \Om^0_L\PS^0_R \subset \Dc^0_0$. Then $[\i\eps\nabla,r]\in [\Dc^1_1,\Dc^{-1}_0] \subset \Dc^0_1$ so that finally $\i\eps[\nabla,\nablab] \equiv -\i\eps\frac{\d}{\d x^i}\frac{\d}{\d p_i} \mod \Dc^0_1$. In conclusion
$$
D^2\equiv -\i\eps \frac{\d}{\d x^i}\frac{\d}{\d p_i} \mod (\Dc^0_1 + \Dc^0_2) \equiv -\i\eps \frac{\d}{\d x^i}\frac{\d}{\d p_i} \mod \Dc^0_1
$$
which shows that $-D^2$ is a generalized Laplacian. \cqfd\\

\begin{example}\label{ecan}\textup{Let $d$ be the exterior derivative of differential forms over $M$. Hence $d\in \PS(M,E)$ is a differential operator of order one. Its right multiplication on $\CS(M,E)$ defines an element of odd degree $d_R\in \Lc(M)$. In a local coordinate system over $U$ one has $d=\i p_i\psi^i$, hence
\be
d_R=\i (p_i\psi^i)_R=\i \psi^i_R p_{iR} = - \psi^i_R\frac{\d}{\d x^i} + \i p_{iL}\psi^i_R
\ee
with $p_{iL}\psi^i_R \in \SPS^1(U,E)_L\Om^1(U)_R$. This shows that $\nabla=- d_R $ is a possible choice. Adding any $\nablab$, the generalized Dirac operator $D=-\i\eps d_R +\nablab$ thus obtained will be called a \emph{de Rham-Dirac} operator on $M$. Note that $\nabla= -d_R$ is completely canonical, only the $\nablab$ part requires some choice.}
\end{example} 

\begin{proposition}
Let $D=-\i\eps d_R +\nablab$ be a de Rham-Dirac operator on $M$. In a local coordinate system over an open set $U\subset M$, the associated generalized Laplacian reads
\beq
-D^2 &=&  \i\eps \Big( \frac{\d}{\d x^i}\frac{\d}{\d p_i} + \sum_{|\al|=2}^{\infty} (a_{\al}^i)_L\frac{\d}{\d x^i} \Big(\frac{\d}{\d p}\Big)^{\al} \Big) \non\\
&& + \eps \Big( p_{iL} \frac{\d}{\d p_i} + \sum_{|\al|=2}^{\infty} (a_{\al}^ip_i)_L \Big(\frac{\d}{\d p}\Big)^{\al} \Big) \label{lap1}\\
 && + \eps \Big((\psi^i\psib_i)_R  + \sum_{|\al|=1}^{\infty} (b_{\al j}^i)_L (\psi^j\psib_i)_R\Big(\frac{\d}{\d p}\Big)^{\al} \Big) \non
\eeq
where $a_{\al}^i, b_{\al j}^i \in \Om^0(U)$ are scalar functions. 
\end{proposition}
{\it Proof:} Since $d^2=0$ and $\nablab^2=0$ one has $-D^2=\i\eps [d_R,\nablab]$. In a local coordinate system one can write $\nablab\equiv \psib_{iR}\frac{\d}{\d p_i} \mod \Om^0_L\Vect_R\cap \Dc^{-1}_0$. Let us calculate
$$
[d_R,\psib_{iR}\frac{\d}{\d p_i}] = \i [(p_j\psi^j)_R , \psib_{iR}\frac{\d}{\d p_i}] = -\i [\psib_i,p_j\psi^j]_R \frac{\d}{\d p_i} - \i \psib_{iR} [(p_j\psi^j)_R  , \frac{\d}{\d p_i} ]\ .
$$
One has $[\psib_i,p_j\psi^j]_R = (p_j[\psib_i,\psi^j])_R = p_{iR}$ and $\psib_{iR} [(p_j\psi^j)_R  , \frac{\d}{\d p_i} ] = - \psib_{iR}\psi^i_R = (\psi^i\psib_i)_R$, so that
$$
[d_R,\psib_{iR}\frac{\d}{\d p_i}] = -\i p_{iR} \frac{\d}{\d p_i} -\i (\psi^i\psib_i)_R = \frac{\d}{\d x^i}\frac{\d}{\d p_i} -\i p_{iL} \frac{\d}{\d p_i} -\i (\psi^i\psib_i)_R\ .
$$
This gives the three principal terms in (\ref{lap1}). The other terms, which are perturbations, come from the commutator of $d_R$ with $\Om^0_L\Vect_R\cap \Dc^{-1}_0$. Indeed, a generic element in $\Om^0_L\Vect_R\cap \Dc^{-1}_0$ can be expanded as $\sum_{|\al|=2}^{\infty} (a_{\al}^i)_L \psib_{iR}\d_p^{\al}$, with $a_{\al}^i \in\Om^0$. One has
$$
[d_R, (a_{\al}^i)_L \psib_{iR}\d_p^{\al}] =  (a_{\al}^i)_L \big( [d_R , \psib_{iR}] \d_p^{\al} - \psib_{iR} [d_R , \d_p^{\al}] \big)
$$
where $[d_R , \psib_{iR}] = \i [(p_j\psi^j)_R , \psib_{iR}] = -\i [\psib_i,\psi^j]_R p_{jR} = -\i p_{iR} =\frac{\d}{\d x^i} -\i p_{iL}$. Moreover $[d_R , \d_p^{\al}]=\i [(p_j\psi^j)_R , \d_p^{\al}] = \i\psi^j_R [p_{jR} , \d_p^{\al}]$ is a sum of terms proportional to $\psi^j_R\d_p^{\beta}$ for all multi-indices $\beta$ such that $|\beta|=|\al|-1$. Hence we can write
$$
[d_R, (a_{\al}^i)_L \psib_{iR}\d_p^{\al}] = (a_{\al}^i)_L\frac{\d}{\d x^i}\d_p^{\al} -\i (a_{\al}^ip_i)_L \d_p^{\al} - \i \sum_{|\beta|=|\al|-1} (b_{\beta j}^i)_L (\psi^j\psib_i)_R\d_p^{\beta}
$$
where the terms of the right-hand-side contribute to the first, second and third line of (\ref{lap1}) respectively. \cqfd\\

\begin{remark}\textup{For any generalized Dirac operator $D=\i\eps\nabla+\nablab$, we can write
\be
\nabla = - d_R + s \quad \mbox{with} \quad s\in \SPS^1(M,E)_L\Om^1(M)_R
\ee
globally on $M$. This property completely characterizes the class of operators $\nabla$ without reference to any local coordinate system. }
\end{remark}

\begin{example}\textup{We now give another important example of generalized Dirac operator related to a choice of torsion-free affine connection $\Gamma$ on $M$. Such a connection is characterized in any local coordinate system over $U\subset M$ by its Christoffel symbols $\Gamma^k_{ij}(x)$, for $i,j,k=1,\ldots,n$, which are symmetric with respect to the lower indices $ij$. Under a coordinate transformation $x^i\mapsto \gamma(x^i)=y^i$ the Christoffel symbols change according to
\be
\Gamma^k_{ij}(x)\mapsto \ ^{\gamma}\Gamma^k_{ij}(x) = \frac{\d x^k}{\d y^l}\frac{\d^2 y^l}{\d x^i\d x^j} + \frac{\d x^k}{\d y^l}\frac{\d y^p}{\d x^i}\frac{\d y^q}{\d x^j} \Gamma^l_{pq}(y)\ .
\ee
In the given coordinate system we define a ``covariant derivative'' operator acting on $\CS(U,E)$:
\be
\nabla_i^{\Gamma} = \frac{\d}{\d x^i} + \big(\Gamma^k_{ij}(x)\big)_{\! L}\Big(p_{kL}\frac{\d}{\d p_j} +  (\psib_k\psi^j)_L - (\psib_k\psi^j)_R\Big) \ .
\ee
Note that it is not quite  a derivation on the algebra $\CS(U,E)$, because $x$ and $p$ do not commute, however its action on the generators $x,p,\psi,\psib$ is what we expect from a covariant derivative:
$$
\nabla_i^{\Gamma}(x^k) = \delta_i^k\ ,\quad  \nabla_i^{\Gamma}(p_j) = \Gamma^k_{ij}p_k\ ,\quad  \nabla_i^{\Gamma}(\psi^k) = -\Gamma^k_{ij}\psi^j\ ,\quad  \nabla_i^{\Gamma} = \Gamma^k_{ij} \psib_k\ .
$$
We say that a generalized Dirac operator $D=\i\eps\nabla+\nablab$ is \emph{affiliated to the connection} $\Gamma$ is in any coordinate system one has
\be
\nabla = \psi^i_R\nabla_i^{\Gamma} + s  \label{affil}\ ,
\ee
where the remainder $s$ has an expansion of the form
\be
s = \psi^i_R\Big(\sum_{|\al|=2}^{\infty} (s^k_{\al i}p_k)_L \d_p^{\al} + \sum_{|\al|=1}^{\infty} (s^k_{\al ij}\psib_k\psi^j + s_{\al i})_L \d_p^{\al} \Big) \label{s}
\ee
for some scalar functions $s^k_{\al i}, s^k_{\al ij}, s_{\al i}\in \Om^0(U)$. Observe that $s$ belongs to $\SPS^1(U,E)_L\Om^1(U)_R\cap \Dc^0_0(U)$. In order to check that this definition makes sense, one has to inspect the transformation law of $\psi^i_R\nabla^{\Gamma}_i$ under a coordinate change $\gamma$. Using the symmetry $\Gamma^k_{ij} = \Gamma^k_{ji}$ one has
$$
\psi^i_R\nabla^{\Gamma}_i = \psi^i_R\frac{\d}{\d x^i} + \psi^i_R\big(\Gamma^k_{ij}(x)p_k\big)_{\! L} \frac{\d}{\d p_j} +  \psi^i_R\big(\Gamma^k_{ij}(x)\psib_k\psi^j\big)_{\! L}  \ .
$$
We already know that $\gamma(\psi^i_R\frac{\d}{\d x^i})\equiv \psi^i_R\frac{\d}{\d x^i}  \mod \SPS^1(U,E)_L\Om^1(U)_R$, but a closer examination of Equation (\ref{change}) 
$$
\gamma\Big( \psi^i_R\frac{\d}{\d x^i} \Big) = \psi^i_R\frac{\d}{\d x^i} -\i p_{iL}\psi^i_R + \Big(\i \frac{\d x^k}{\d y^l}p_k + \frac{\d}{\d x^q}\frac{\d x^k}{\d y^l}\psi^q\psib_k\Big)_{\!\! L} \Big(\frac{\d y^l}{\d x^i}\psi^i\Big)_{\!\! R}
$$
gives, by means of the expansion $\big(\frac{\d y^l}{\d x^i}\psi^i\big)_{\! R} = \sum_{|\al|=0}^{\infty} \frac{(-\i)^{|\al|}}{\al !} \d_x^{\al}\big(\frac{\d y^l}{\d x^i}\big)_{\! L}\psi^i_R\d_p^{\al}$, 
\beq
\gamma\Big( \psi^i_R\frac{\d}{\d x^i} \Big) &=& \psi^i_R\frac{\d}{\d x^i} + \Big(\frac{\d x^k}{\d y^l}p_k\frac{\d^2 y^l}{\d x^i\d x^j}\Big)_{\!\! L}\psi^i_R\frac{\d}{\d p_j} + \Big(\frac{\d x^k}{\d y^l}\frac{\d^2 y^l}{\d x^i\d x^j}\psib_k\psi^j\Big)_{\!\! L}\psi^i_R \non\\
&& + \sum_{|\al|=2}^{\infty} \frac{(-\i)^{|\al|}}{\al !} \Big(\i\frac{\d x^k}{\d y^l}p_k\d_x^{\al}\frac{\d y^l}{\d x^i}\Big)_{\!\! L}\psi^i_R\d_p^{\al} \non\\
&& + \sum_{|\al|=1}^{\infty} \frac{(-\i)^{|\al|}}{\al !} \Big(\frac{\d}{\d x^q}\frac{\d x^k}{\d y^l}\psi^q\psib_k\d_x^{\al}\frac{\d y^l}{\d x^i}\Big)_{\!\! L}\psi^i_R\d_p^{\al} \ .\non
\eeq
We used the identities $-\i p_i +\i \frac{\d x^k}{\d y^l}p_k \frac{\d y^l}{\d x^i} = \frac{\d x^k}{\d y^l} \frac{\d^2 y^l}{\d x^k\d x^i}$ and $\psi^j\psib_k = \delta^j_k-\psib_k\psi^j$ in order to simplify the first line. Since commutators with $p$ are proportional to derivations with respect to $x$, the above expression reads
$$
\gamma\Big( \psi^i_R\frac{\d}{\d x^i} \Big) = \psi^i_R\frac{\d}{\d x^i} + \psi^i_R\Big(\frac{\d x^k}{\d y^l}\frac{\d^2 y^l}{\d x^i\d x^j}\Big)_{\!\! L} \Big(p_{kL}\frac{\d}{\d p_j} + (\psib_k\psi^j)_L\Big) + s' \ ,
$$
where the remainder $s'$ has an expansion of the form (\ref{s}). In the same way, one can show that
\beq
\lefteqn{\gamma\Big(\psi^i_R\big(\Gamma^k_{ij}(x)\big)_{\! L} \Big(p_{kL}\frac{\d}{\d p_j} + (\psib_k\psi^j)_L\Big)\Big) =} \non\\
&& \qquad \psi^i_R\Big(\frac{\d x^k}{\d y^l}\frac{\d y^p}{\d x^i}\frac{\d y^q}{\d x^j} \Gamma^l_{pq}(y)\Big)_{\!\! L} \Big(p_{kL}\frac{\d}{\d p_j} + (\psib_k\psi^j)_L\Big) + s'' \non
\eeq
with a remainder $s''$ of the form (\ref{s}). Hence $\gamma(\psi^i_R\nabla_i^{\Gamma}) = \psi^i_R\nabla_i^{^{\gamma}\Gamma} + s$, and using a partition of unity we can build a global operator $\nabla$ on $M$ with the wanted property. The following proposition, which is an analogue of the Lichnerowicz formula, relates the square of the corresponding Dirac operator to the curvature tensor of the connection $\Gamma$, whose components in local coordinates are }
\be
R^k_{lij} = \frac{\d\Gamma^k_{jl}}{\d x^i}  - \frac{\d\Gamma^k_{il}}{\d x^j}  + \Gamma^k_{im}\Gamma^m_{jl} - \Gamma^k_{jm}\Gamma^m_{il}\ .
\ee
\end{example}

\begin{proposition}
Let $D$ be a Dirac operator affiliated to a torsion-free affine connection $\Gamma$ on $M$. In a local coordinate system over an open set $U\subset M$, the associated generalized Laplacian reads
\beq
-D^2 &=&  \i\eps \Big( \frac{\d}{\d x^i}\frac{\d}{\d p_i} + (\Gamma^k_{ij})_L(\psi^i\psib_k)_R \frac{\d}{\d p_j} + u + v \Big) \non\\
&& + \eps^2 \Big( \frac{1}{2}\, (\psi^i\psi^j)_R (R^k_{lij})_L  \Big( p_{kL} \frac{\d}{\d p_l} + (\psib_k\psi^l)_L\Big) + w \Big) \label{lap2}
\eeq
where $R^k_{lij}$ are the components of the curvature tensor, and
\beq
u &=&  \sum_{|\al|=2}^{\infty} \Big( (u_{\al i})_L\frac{\d}{\d x^i} + (u^k_{\al}p_k)_L + (u^k_{\al i})_L (\psi^i\psib_k)_R + (u_{\al})_L \Big) \d_p^{\al} \non\\
v &=& \sum_{|\al|=1}^{\infty} (v^k_{\al i} \psib_k\psi^i)_L \d_p^{\al} \non\\
w &=& (\psi^i\psi^j)_R\Big(\sum_{|\al|=2}^{\infty} (w^k_{\al ij}p_k)_L \d_p^{\al} + \sum_{|\al|=1}^{\infty} (w^k_{\al lij}\psib_k\psi^l + w_{\al ij})_L \d_p^{\al} \Big) \non
\eeq
where $u_{\al i}, u^k_{\al}, u^k_{\al i}, u_{\al}, v^k_{\al i}, w^k_{\al ij}, w^k_{\al lij}, w_{\al ij} \ni\Om^0(U)$ are scalar functions.
\end{proposition}
{\it Proof:} Since $\nablab^2=0$ one has $-D^2=-\i\eps[\nabla,\nablab] + \eps^2\nabla^2$. In a local coordinate system $\nabla = \psi^i_R\nabla_i^{\Gamma} + s$ and $\nablab = \psib_{kR}\frac{\d}{\d p_k} + r$ with 
\beq
s &=& \psi^i_R\Big(\sum_{|\al|=2}^{\infty} (s^k_{\al i}p_k)_L \d_p^{\al} + \sum_{|\al|=1}^{\infty} (s^k_{\al ij}\psib_k\psi^j + s_{\al i})_L \d_p^{\al} \Big) \non\\
r &=& \sum_{|\al|=2}^{\infty} (r^i_{\al})_L \psib_{iR} \d_p^{\al}\ . \non
\eeq
Hence $[\nabla,\nablab] = [\psi^i_R\nabla_i^{\Gamma} , \psib_{kR}\frac{\d}{\d p_k}] + [\psi^i_R\nabla_i^{\Gamma} , r] + [s, \psib_{kR}\frac{\d}{\d p_k}] + [s,r]$. We compute each commutator of the right hand side separately. Firstly,
\beq
\lefteqn{[\psi^i_R\nabla_i^{\Gamma} , \psib_{kR}\frac{\d}{\d p_k}] =} \non\\
&& [ \psi^i_R \frac{\d}{\d x^i}, \psib_{kR}\frac{\d}{\d p_k} ] + [\big(\Gamma^l_{ij}p_l\big)_{\! L}\psi^i_R\frac{\d}{\d p_j} , \psib_{kR}\frac{\d}{\d p_k} ] + [\big(\Gamma^l_{ij}\psib_l\psi^j\big)_{\! L}\psi^i_R , \psib_{kR}\frac{\d}{\d p_k} ] \ . \non
\eeq
One has 
$$
[ \psi^i_R \frac{\d}{\d x^i} , \psib_{kR}\frac{\d}{\d p_k} ] = [ \psi^i_R , \psib_{kR}] \frac{\d}{\d x^i}\frac{\d}{\d p_k} = - [\psib_k , \psi^i]_R \frac{\d}{\d x^i}\frac{\d}{\d p_k} = - \frac{\d}{\d x^i}\frac{\d}{\d p_i}\ .
$$ 
Then
\beq
\lefteqn{[\big(\Gamma^l_{ij}p_l\big)_{\! L}\psi^i_R\frac{\d}{\d p_j} , \psib_{kR}\frac{\d}{\d p_k} ] } \non\\
&=& \big(\Gamma^l_{ij}p_l\big)_{\! L} [\psi^i_R,\psib_{kR}] \frac{\d}{\d p_k}\frac{\d}{\d p_j} - \psib_{kR}\big(\Gamma^l_{ij}\big)_{\! L} [ p_{lL} ,\frac{\d}{\d p_k} ] \psi^i_R\frac{\d}{\d p_j} \non\\
&=& - \big(\Gamma^l_{ij}p_l\big)_{\! L} \frac{\d}{\d p_i}\frac{\d}{\d p_j} - \big(\Gamma^k_{ij}\big)_{\! L} (\psi^i\psib_k)_R \frac{\d}{\d p_j} \non
\eeq
and
$$
 [\big(\Gamma^l_{ij}\psib_l\psi^j\big)_{\! L}\psi^i_R , \psib_{kR}\frac{\d}{\d p_k} ] = - \big(\Gamma^l_{ij}\psib_l\psi^j\big)_{\! L} \frac{\d}{\d p_i}
$$
so that
\beq
-\i\eps [ \psi^i_R \nabla_i^{\Gamma} , \psib_{kR}\frac{\d}{\d p_k} ] &=& \i\eps \frac{\d}{\d x^i}\frac{\d}{\d p_i} + \i\eps (\Gamma^k_{ij})_L (\psi^i\psib_k)_R \frac{\d}{\d p_j} \non\\
&&  +\i\eps \big(\Gamma^l_{ij}p_l\big)_{\! L} \frac{\d}{\d p_i}\frac{\d}{\d p_j} +\i\eps \big(\Gamma^l_{ij}\psib_l\psi^j\big)_{\! L} \frac{\d}{\d p_i} \ . \non
\eeq
The first and second term appear in the first line of (\ref{lap2}), while the third and fourth terms contribute to $u$ and $v$ respectively. We continue with the commutator $[\psi^i_R\nabla_i^{\Gamma} , r]$:
$$
[\psi^i_R \frac{\d}{\d x^i}, r] = \sum_{|\al|=2}^{\infty} [\psi^i_R \frac{\d}{\d x^i}, (r^j_{\al})_L \psib_{jR} ] \d_p^{\al} = \sum_{|\al|=2}^{\infty} \Big(\Big(\frac{\d r^j_{\al}}{\d x^i}\Big)_{\!\! L} \psi^i_R\psib_{jR} - (r^i_{\al})_L\frac{\d}{\d x^i}\Big) \d_p^{\al}
$$
and
\beq
[\big(\Gamma^l_{ij}p_l\big)_{\! L}\psi^i_R\frac{\d}{\d p_j} , r ] &=& \sum_{|\al|=2}^{\infty} [\big(\Gamma^l_{ij}p_l\big)_{\! L}\psi^i_R , (r^j_{\al})_L \psib_{jR} \d_p^{\al}]\frac{\d}{\d p_j}  \non\\
&=& \sum_{|\al|=3}^{\infty} (a^k_{\al}p_k)_L\d_p^{\al} + \sum_{|\al|=2}^{\infty} \big((a^k_{\al i})_L (\psi^i\psib_k)_R + (a_{\al} )_L\big) \d_p^{\al}   \non
\eeq
for some scalar functions $a^k_{\al}, a^k_{\al i}, a_{\al}$, and
$$
[\big(\Gamma^l_{ij}\psib_l\psi^j\big)_{\! L}\psi^i_R , r] = - \sum_{|\al|=2}^{\infty} \big(\Gamma^l_{ij}\psib_l\psi^j r^i_{\al}\big)_{\! L} \d_p^{\al}\ .
$$
Hence $[\psi^i_R\nabla_i^{\Gamma} , r]$ can be absorbed inside $u+v$. Further on, we have
\beq 
[\psib_{jR}\frac{\d}{\d p_j}, s] &=&  \sum_{|\al|=2}^{\infty} [\psib_{jR}\frac{\d}{\d p_j}, (s^k_{\al i}p_k)_L \psi^i_R ] \d_p^{\al}  \non\\
&+& \sum_{|\al|=1}^{\infty} (s^k_{\al il}\psib_k\psi^l)_L [\psib_{jR}\frac{\d}{\d p_j}, \psi^i_R] \d_p^{\al} + \sum_{|\al|=1}^{\infty} (s_{\al i})_L [\psib_{jR}\frac{\d}{\d p_j}, \psi^i_R ] \d_p^{\al}   \non\\
&=&   \sum_{|\al|=2}^{\infty} \Big( (s^j_{\al i})_L \psib_{jR}\psi^i_R  - (s^k_{\al i}p_k)_L \frac{\d}{\d p_i} \Big) \d_p^{\al} \non\\
&& - \sum_{|\al|=1}^{\infty} (s^k_{\al il}\psib_k\psi^l)_L \frac{\d}{\d p_i} \d_p^{\al} - \sum_{|\al|=1}^{\infty} (s_{\al i})_L \frac{\d}{\d p_i} \d_p^{\al}   \non
\eeq
The first and third series of the right-hand-side can be absorbed inside $u$, whereas the second series counts for $v$. Instead of computing the commutator $[s,r]$ explicitly, we only need to remark that $s\in \SPS^1_L\Om^1_R\cap \Dc^0_0$ and $r\in \Om^0_L\Vect_R\cap \Dc^{-1}_0$. Then
\beq
{[\SPS^1_L\Om^1_R, \Om^0_L\Vect_R ]} &\subset& [\SPS^1_L, \Om^0_L]\PS^0_R + \SPS^1_L[\Om^1_R,\Vect_R ] \non\\
&\subset& \Om^0_L\PS^0_R + \SPS^1_L\Om^0_R  \non
\eeq
It follows that $[s,r] \in (\Om^0_L\PS^0_R + \SPS^1_L\Om^0_R)\cap\Dc^{-1}_0$ can be absorbed inside $u+v$. Now we look at 
$$
\nabla^2 = (\psi^i_R\nabla_i^{\Gamma}+s)^2 = (\psi^i_R\nabla_i^{\Gamma})^2 + [\psi^i_R\nabla_i^{\Gamma}, s] + s^2\ .
$$
A routine computation gives
$$
[\nabla_i^{\Gamma}, \nabla_j^{\Gamma}] = (R^k_{lij})_L \Big( p_{kL}\frac{\d}{\d p_l} + (\psib_k\psi^l)_L - (\psib_k\psi^l)_R \Big)\ .
$$ 
Consequently, the Bianchi identity $(R^k_{lij})_L (\psi^l\psi^i\psi^j)_R=0$ implies 
$$
(\psi^i_R\nabla_i^{\Gamma})^2 =  \frac{1}{2}\,\psi^i_R\psi^j_R [\nabla_i^{\Gamma}, \nabla_j^{\Gamma}] = \frac{1}{2}\, (\psi^i\psi^j)_R (R^k_{lij})_L  \Big( p_{kL} \frac{\d}{\d p_l} + (\psib_k\psi^l)_L\Big) \ .
$$
This the leading term in the second line of (\ref{lap2}). Then we have 
$$
[\psi^l_R \frac{\d}{\d x^l}, s] = (\psi^l\psi^i)_R\Big(\sum_{|\al|=2}^{\infty} \Big(\frac{\d s^k_{\al i}}{\d x^l}p_k\Big)_{\!\! L} \d_p^{\al} + \sum_{|\al|=1}^{\infty} \Big(\frac{\d s^k_{\al ij}}{\d x^l}\psib_k\psi^j + \frac{\d s_{\al i}}{\d x^l}\Big)_{\!\! L} \d_p^{\al} \Big) 
$$
and 
$$
[\big(\Gamma^l_{ij}p_l\big)_{\! L}\psi^i_R\frac{\d}{\d p_j} , s ] = (\psi^i\psi^j)_R\Big(\sum_{|\al|=2}^{\infty} (b^k_{\al ij}p_k)_L \d_p^{\al} + \sum_{|\al|=1}^{\infty} (b^k_{\al lij}\psib_k\psi^l + b_{\al ij})_L \d_p^{\al} \Big) 
$$
for some scalar functions $b^k_{\al ij}, b^k_{\al lij}, b_{\al ij}$, and 
$$
[\big(\Gamma^l_{ij}\psib_l\psi^j\big)_{\! L}\psi^i_R , s ] = (\psi^i\psi^j)_R \sum_{|\al|=1}^{\infty} (c^k_{\al lij}\psib_k\psi^l)_L \d_p^{\al}
$$
for some other scalar functions $c^k_{\al lij}$. Hence $[\psi^i_R\nabla_i^{\Gamma}, s]$ can be absorbed inside $w$. Finally one easily checks that $s^2$ is also of the form $w$. \cqfd\\

\section{Algebraic JLO formula}\label{sjlo}

We first recall Connes' definition of periodic cyclic cohomology \cite{C86}. Let $\Ac$ be a trivially-graded associative $\cc$-algebra. Form the unitalized algebra $\Ac^+=\Ac\oplus\cc$, even if $\Ac$ already has a unit. For any $k\in\nn^*$ denote by $CC^k(\Ac)$ the space of $(k+1)$-linear maps $\Ac^+\times\Ac^{\times k}\to\cc$, and $CC^0(\Ac)$ the space of linear maps $\Ac\to\cc$. The Hochschild operator $b:CC^k(\Ac)\to CC^{k+1}(\Ac)$ is defined on a $k$-cochain $\varphi_k\in CC^k(\Ac)$ by
\beq
b\varphi_k(a_0,\ldots,a_{k+1}) &=& \sum_{i=0}^k (-1)^i \varphi_k(a_0,\ldots,a_ia_{i+1},\ldots,a_{k+1})\non\\
&& + (-1)^{k+1} \varphi_k(a_{k+1}a_0,\ldots,a_{k})
\eeq
for any $a_0\in\Ac^+$ and $a_1,\ldots,a_k\in \Ac$. The Connes operator $B:CC^k(\Ac)\to CC^{k-1}(\Ac)$ reads
\be
B\varphi_k(a_0,\ldots,a_{k-1}) = \sum_{i=0}^{k-1} (-1)^{i(k-i)} \varphi_k(a_i,\ldots,a_{k-1}, a_0,\ldots,a_{i-1})\ .
\ee
One checks $b^2=B^2=bB+Bb=0$. The direct sum $CP^{\bullet}(\Ac)=\sum_{k=0}^{\infty}CC^k(\Ac)$ endowed with the boundary operator $b+B$ is therefore a $\zz_2$-graded complex. The cohomology $HP^{\bullet}(\Ac)$, of this complex is the periodic cyclic cohomology of $\Ac$. Thus, an even periodic cyclic cocycle over $\Ac$ is a \emph{finite} collection $\varphi=(\varphi_0,\varphi_2, \ldots,\varphi_{2n})$ of homogeneous cochains such that
\be
b\varphi_k+B\varphi_{k+2}=0\quad \mbox{for $0\leq k<2n$}\ ,\quad b\varphi_{2n}=0\ .
\ee
An odd periodic cyclic cocycle is a finite collection $\varphi=(\varphi_1,\varphi_3, \ldots,\varphi_{2n+1})$ verifying analogous relations.

\begin{example}\textup{(Connes \cite{C86}) If $M$ is a compact manifold, any homology class $[C_k]\in H_k(M,\cc)$ represented by a $k$-dimensional closed de Rham current $C_k$ gives rise to a periodic cyclic cohomology class over the commutative algebra $\cinf(M)$ by setting
\be
\varphi_k(a_0,\ldots,a_k)= \frac{c_k}{k!}\, \langle {C_k}, a_0da_1\ldots da_k \rangle \ ,\quad \forall a_i\in\cinf(M)\ ,
\ee
where $c_k$ is a normalization factor depending on the parity of $k$. We choose $c_{2k}= 1/(2\pi\i)^k$ and $c_{2k+1}= 1/(2\pi\i)^{k+1}$ for compatibility with the usual normalization of characteristic classes in de Rham cohomology. Then one checks $b\varphi_k=0=B\varphi_k$ so that $[\varphi_k]\in HP^{k\mod 2}(\cinf(M))$ is represented by a homogeneous cochain of degree $k$. One thus gets a linear map
\be
H_{\bullet}(M,\cc)\to HP^{\bullet}(\cinf(M))\label{derham}
\ee
for any compact manifold. In fact, Connes shows that this is an \emph{isomorphism} \cite{C86}, provided that cyclic cohomology is defined through continuous cochains with respect to the natural locally convex topology of $\cinf(M)$. Since we are not concerned with analytical issues in this paper, the fact that (\ref{derham}) is an isomorphism will be irrelevant for us. }
\end{example}

\begin{example}\textup{Consider the non-commutative algebra $\CS^0(M)$ of formal symbols of order $\leq 0$ on a closed manifold $M$. The leading symbol gives rise to an algebra homomorphism $\la:\CS^0(M)\to \cinf(S^*M)$ to the commutative algebra of functions over the cosphere bundle $S^*M$. Since cyclic cohomology pullbacks under homomorphisms, one gets, modulo composition with (\ref{derham}), a canonical map
\be
\la^*: H_{\bullet}(S^*M)\to HP^{\bullet}(\CS^0(M))\ . \label{iso}
\ee
In fact, Wodzicki shows that this is an \emph{isomorphism} \cite{W88}, provided the natural locally convex topology of $\CS^0(M))$ is taken into account. Again, we will not use the fact that $\la^*$ is an isomorphism. }
\end{example}

Now fix a closed $n$-dimensional manifold $M$. We will construct some cyclic cocycles over the algebra $\CS^0(M)$ using Dirac operators as defined in section \ref{ddirac}. By construction $\CL^0(M)$ is an algebra of operators on the space $\cinf(M)$. We can view $\CL^0(M)$ as an algebra of operators on the space of sections of the vector bundle $E=\Lambda T^*_{\cc}M$: indeed its action on the zero-forms $\cinf(M)=\Om^0(M)$ can be extended by zero on $\Om^k(M)$, $\forall k\geq 1$. Therefore one has a canonical homomorphism of $\CL^0(M)$ into the even part of the $\zz_2$-graded algebra $\CL^0(M,E)$. It descends to an homomorphism $\pi: \CS^0(M)\to \CS^0(M,E)$. In a local coordinate system we can write
\be
\pi(a)(x,p,\psi,\psib) = a(x,p) \Pi \qquad \forall a\in \CS^0(M)\ ,
\ee
where $\Pi = \psib_1\psi^1\ldots \psib_n\psi^n$ is the Clifford section corresponding to the projection operator from $\Om^*(M)$ onto $\Om^0(M)$. Then we can compose $\pi$ with the left representation of $\CS^0(M,E)$ as endomorphisms on the vector space $\CS(M,E)$. This yields an injective homomorphism of algebras 
\be
\rho: \CS^0(M)\hookrightarrow \Dc^0_0(M)\ ,\quad \rho(a)= (a\Pi)_L\quad \forall a\in \CS^0(M)\ .
\ee
We are now ready to introduce the following algebraic version of the JLO cocycle \cite{JLO}. It involves the graded trace on the algebra of trace-class operators $\Tc(M)$ introduced in section \ref{sct}.
\begin{proposition}\label{pJLO1}
Let $D=\i\eps\nabla+\nablab \in \Dc^1(M)$ be a generalized Dirac operator. The homogeneous cochains over the algebra $\CS^0(M)$
\be
\varphi_k^D(a_0,\ldots,a_k)= \int_{\Delta_k} \Tr_s \big(\rho(a_0)e^{-t_0D^2}[D,\rho(a_1)] e^{-t_1D^2} \ldots [D,\rho(a_k)]e^{-t_kD^2}\big) dt
\ee
defined for all $k\in 2\nn$, are the components of an even periodic cyclic cocycle $\varphi^D$ and vanish whenever $k>2n$, $n=\dim M$. Moreover, the periodic cyclic cohomology class $[\varphi^D]\in HP^0(\CS^0(M))$ does not depend on $D$. 
\end{proposition}
{\it Proof:} The graded trace of a trace-class operator $s\in \Tc(M)$ vanishes if the Clifford part of $s$ is not of heighest weight, that is, if $s$ is not proportional to the product $(\psi^1\ldots\psi^n\psib_1\ldots\psib_n)_L (\psi^1\ldots\psi^n\psib_1\ldots\psib_n)_R$ in local coordinates. Hence in the computation of $\varphi^D_k$, we should only retain the terms which bring at least $n$ powers of $\psi_L$ (resp. of $\psi_R$) and exactly the same powers of $\psib_L$ (resp. of $\psib_R$), because we have to take into account the possible lowering of powers coming from commutators $[\psi^i,\psib_j]=\delta^i_j$. All other combinations of $\psi_L,\psib_L,\psi_R,\psib_R$ will vanish under the graded trace. In fact the right sector $\psi_R,\psib_R$ will be our main interest. One has
$$
[D,\rho(a)] = \i\eps[\nabla, (a\Pi)_L] + [\nablab,(a\Pi)_L]\ .
$$
The first term brings a factor $\eps\psi_R$, whereas the second term brings a factor $\psib_R$. We define the \emph{pseudodifferential order} of an operator according to the following rule: $a_L$ is of order $m$ for any symbol $a\in\CS^m(M,E)$, the operators $\psi_R,\psib_R, \eps$ are of order $0$, while $\d_p$ is of order $-1$ and $\d_x$ of order $+1$. From these rules one sees that the operator $\i\eps[\nabla, (a\Pi)_L]$ has order $\leq 0$, and $[\nablab,(a\Pi)_L]$ has order $\leq -1$. In the same way we inspect the generalized Laplacian
$$
-D^2 = -\i\eps [\nabla,\nablab] + \eps^2 \nabla^2 \ .$$
From the proof of Proposition \ref{plap} we know that $\nabla^2\in \SPS^1_L\Om^2_R$, hence $\eps^2\nabla^2$ has pseudodifferential order $\leq 1$ and brings a factor $\eps^2\psi_R\psi_R$. Similarly one has $-\i\eps[\nabla,\nablab]=\Delta + u$ where $\Delta=\i\eps\frac{\d}{\d x^i}\frac{\d}{\d p_i}$ is the flat Laplacian in local coordinates. $u$ has order $\leq 0$ and its right sector is proportional to either $\eps\psi_R\psib_R$ or $1$. We treat $-D^2$ as a perturbation of the flat Laplacian. A Duhamel expansion of the exponentials $\exp(-t_iD^2)$ appearing in the cochain $\varphi^D$ leads to the computation of terms like
\beq
\lefteqn{\Tr_s\big( \rho(a_0) \exp(t_0\Delta) \, X_1 \, \exp(t_1\Delta)\ldots X_k \, \exp(t_k\Delta) \big) =} \non\\
&&\qquad \qquad  \bint \big\langle\!\big\langle (a_0\Pi)_L \si^{t_0}_{\Delta}(X_1) \ldots \si^{t_0+\ldots +t_{k-1}}_{\Delta}(X_k) \, \exp \Delta  \big\rangle\!\big\rangle [n] \non
\eeq
where $X_i=\eps^2\nabla^2$, or $X_i=u$, or $X_i=\i\eps[\nabla,(a_j\Pi)_L]$, or $X_i=[\nablab,(a_j\Pi)_L]$ for some $a_j\in \CS^0(M)$. In order to achieve an exact balance between the powers of $\psi_R$ and $\psib_R$, we see that the number $\bar{l}$ of factors $[\nablab,(a_j\Pi)_L]$ should equal $l+2m$, where $l$ is the number of factors $\i\eps[\nabla,(a_j\Pi)_L]$ and $m$ the number of factors $\eps^2\nabla^2$.
The pseudodifferential order of each $X_i$ is not modified by the action of the modular group $\si_{\Delta}$ because
$$
[\Delta, X_i] = \i\eps \Big( \frac{\d X_i}{\d x^j}\frac{\d}{\d p_j} + \frac{\d X_i}{\d p_j}\frac{\d}{\d x^j} + \frac{\d^2 X_i}{\d x^j\d p_j} \Big)\ .
$$
The contractions $\langle\d_x^{\al}\d_p^{\beta}\exp\Delta\rangle$ also preserve the pseudodifferential order ($\d_x$ and $\d_p$ are simultaneously contracted). It follows that the pseudodifferential order of the symbol $\big\langle\!\big\langle \rho(a_0) \si^{t_0}_{\Delta}(X_1) \ldots \si^{t_0+\ldots +t_{k-1}}_{\Delta}(X_k) \, \exp \Delta  \big\rangle\!\big\rangle [n]$ is $\leq -\bar{l} + m = -l-m$, and its Wodzicki residue vanishes unless $-l-m\geq -n$ ($n=\dim M$). The latter condition implies $l \leq n-m$ and $\bar{l} \leq n+m$, so $l+\bar{l}\leq 2n$. This means that $\varphi_k^D$ vanishes whenever it involves more than $2n$ commutators $[D,\rho(a)]$, that is, whenever $k>2n$.\\
Hence $\varphi^D$ is a cochain in the periodic complex $CP^{\bullet}(\CS^0(M))$. The cocycle identity $b\varphi_k^D+B\varphi_{k+2}^D=0$ then follows from well-known algebraic manipulations which we do not need to reproduce here, see \cite{JLO}. Finally observe that given two operators $D_0$ and $D_1$ the linear homotopy
$$
D=tD_1 + (1-t)D_0\ ,\qquad t\in[0,1]\ ,
$$
is a Dirac operator for all $t$. It is again a classical result that the cocycles $\varphi^{D_0}$ and $\varphi^{D_1}$ are related by a transgression formula of JLO type (see for instance \cite{EFJL}). One shows as above that the transgressed cochain, in our case, lies in the periodic complex. Hence the periodic cyclic cohomology class of $\varphi^D$ does not depend on $D$. \cqfd\\

\begin{proposition}\label{pwod}
Let $D=-\i\eps d_R + \nablab$ be a de Rham-Dirac operator. Then $\varphi^D_0$ is the Wodzicki residue on $\CS^0(M)$, while the other components $\varphi^D_k$ vanish for $k>0$. Hence $[\varphi^D]$ is the periodic cyclic cohomology class of the Wodzicki residue. 
\end{proposition}
{\it Proof:} Let us first look at the commutator $[D,\rho(a)]$. Since $\rho(a)=(a\Pi)_L$ belongs to the left sector, it commutes with $d_R$, so that
$$
[D, \rho(a)] = [\nablab, (a\Pi)_L]\ .
$$
By definition $\nablab\in \Om^0(M)_L\Vect(M)_R$ is proportional to $\psib_R$ and not to $\psi_R$. Thus $[D, \rho(a)]$ brings a factor $\psib_R$. On the other hand, the generalized Laplacian $-D^2$ is given by Formula (\ref{lap1}), and brings either $(\psi\psib)_R$ or 1 in the right sector. This means that whenever some commutators $[D,\rho(a)]$ appear, the graded trace must vanish because the $\psib_R$'s cannot be balanced with the same amount of $\psi_R$'s. Hence $\varphi^D_k=0$ whenever $k>0$, and the only remaining component is
$$
\varphi_0^D(a) = \Tr_s( \rho(a) \exp(-D^2))=\Tr_s( (a\Pi)_L \exp(-D^2))\ .
$$
We work in local coordinates $(x,p)$ over $U\subset M$ and suppose that the symbol $a$ has $x$-support contained in $U$ (the general case follows by linearity). Write $-D^2=\Delta + s$, where $\Delta=\i\eps \frac{\d}{\d x^i} \frac{\d}{\d p_i}$ is the canonical flat Laplacian, and the remainder $s$ is given by Equation (\ref{lap1}):
$$
s = \eps \Big( p_{iL} \frac{\d}{\d p_i} + (\psi^i\psib_i)_R + \sum_{|\al|=2}^{\infty} \Big( \i(a_{\al}^i)_L\frac{\d}{\d x^i} + (a_{\al}^ip_i)_L \Big) \d_p^{\al} + \sum_{|\al|=1}^{\infty} (b_{\al j}^i)_L (\psi^j\psib_i)_R \d_p^{\al} \Big)
$$
for some scalar functions $a_{\al}^i, b_{\al j}^i \in \Om^0(U)$. Our goal is to show that the series over the multi-index $\al$ do not contribute to $\varphi^D_0$. We use a Duhamel expansion for $\exp(-D^2)$:
$$
\varphi_0^D(a) =  \sum_{k=0}^{\infty} \int_{\Delta_k} \Tr_s\big( (a\Pi)_L \si^{t_0}_{\Delta}(s)  \si^{t_0+t_1}_{\Delta}(s) \ldots \si^{t_0+\ldots +t_{k-1}}_{\Delta}(s)\,\exp\Delta \big)  dt
$$
Now rewrite the product $\si^{t_0}_{\Delta}(s) \ldots \si^{t_0+\ldots +t_{k-1}}_{\Delta}(s)$ by moving all the derivation operators $\d_x$ and $\d_p$ to the right, in front of $\exp\Delta$. The graded trace would vanish if the resulting powers of $\d_x$ and $\d_p$ are not exactly equal, because it involves the contractions $\langle \d_x \d_p \exp\Delta\rangle$. We remark that all the terms in $s$ except $(\psi^i\psib_i)_R$ bring a power of $\d_p$ strictky higher than the power of $\d_x$. However, a $\d_p$ can be absorbed by commutation with $p_L$ when it moves to the right, and a $\d_x$ can appear from $\si_{\Delta}^t(p_L) = p_L + t[\Delta, p_L] = p_L +\i\eps t\d_x$. A rapid inspection shows that an exact balance between $\d_x$ and $\d_p$ cannot occur if either $( \i(a_{\al}^i)_L\frac{\d}{\d x^i} + (a_{\al}^ip_i)_L ) \d_p^{\al}$ with $|\al|\geq 2$, or $(b_{\al j}^i)_L (\psi^j\psib_i)_R \d_p^{\al}$ with $|\al|\geq 1$ appears. Thus we can keep the only relevant part $\eps( p_{iL} \frac{\d}{\d p_i} + (\psi^i\psib_i)_R) $ of $s$ in the product $\si^{t_0}_{\Delta}(s) \ldots \si^{t_0+\ldots +t_{k-1}}_{\Delta}(s)$, and write
$$
\varphi_0^D(a) = \Tr_s \Big( (a\Pi)_L \exp\big(\Delta + \eps p_L \cdot \d_p + \eps (\psi^i\psib_i)_R\big) \Big)\ .
$$
$(\psi^i\psib_i)_R$ commutes with $\Delta + \eps p_L \cdot \d_p$, hence the exponential splits as the product of $\exp(\eps (\psi^i\psib_i)_R)$ and $\exp(\Delta + \eps p_L \cdot \d_p)$. Expanding $\exp(\eps (\psi^i\psib_i)_R)$ in powers of $\eps$, only the term of order $n$ survives because it involves the product of all $\psi_R$'s and $\psib_R$'s, and the higher powers of $\eps$ are ignored by the graded trace. One finds
\beq
\varphi_0^D(a) &=& \Tr_s \Big( (a\Pi)_L \eps^n (\psi^1\psib_1\ldots \psi^n\psib_n)_R\exp\big(\Delta + \eps p_L \cdot \d_p \big) \Big) \non\\
&=& \bint \tr_s(a\Pi) \big\langle\!\big\langle \eps^n (\psi^1\psib_1\ldots \psi^n\psib_n)_R \exp\big(\Delta + \eps p_L \cdot \d_p \big) \big\rangle\!\big\rangle [n]\ . \non
\eeq
By definition of the graded trace on the Clifford algebra, $\tr_s(a\Pi)=a$ and $\langle (\psi^1\psib_1\ldots \psi^n\psib_n)_R \rangle = (-1)^n \tr_s(\psi^1\psib_1\ldots \psi^n\psib_n)=1$ so that 
$$
\varphi_0^D(a) = \bint a \,  \big\langle \exp\big(\Delta + \eps p_L \cdot \d_p \big)  \big\rangle [0] \ .
$$
Then we apply Lemma \ref{luse} to the matrix $R=\eps\Id$. This yields the formal power series in $\eps$
$$
\langle \exp\big(\Delta + \eps p_L \cdot \d_p \big) \exp(-\Delta) \rangle = \Td(\eps\Id) = \Big( \frac{\eps}{e^{\eps}-1} \Big)^n\ ,
$$ 
whose coefficient of degree zero is  $\Td(\eps\Id)[0]=1$. Therefore $\varphi_0^D(a)$ is the Wodzicki residue as claimed. \cqfd\\

\begin{theorem}\label{twod}
The periodic cyclic cohomology class of the Wodzicki residue vanishes in $HP^0(\CS^0(M))$ for any closed manifold $M$. 
\end{theorem}
{\it Proof:} Let $\Gamma$ be the Levi-Civita connection associated to a given Riemannian metric on $M$, and let $D=\i\eps\nabla +\nablab$ be a generalized Dirac operator affiliated to $\Gamma$. We will show that all the components of the corresponding cocycle $\varphi^D$ vanish. The theorem is then a consequence of Propositions \ref{pJLO1} and \ref{pwod}.\\
In a local coordinate system $\nabla$ is expressed in terms of the Christoffel symbols $\Gamma^k_{ij}$ of the connection:
$$
\nabla = \psi^i_R \frac{\d}{\d x^i} + \big(\Gamma^k_{ij}(x) p_k\big)_L \psi^i_R \frac{\d}{\d p_j} + \big(\Gamma^k_{ij}(x) \psib_k\psi^j\big)_L \psi^i_R + s \ .
$$
The remainder $s$ can be expanded in power series of the partial derivative $\d_p$,
$$
s = \sum_{|\al|=2}^{\infty} (s^k_{\al i}p_k)_L\psi^i_R \d_p^{\al} + \sum_{|\al|=1}^{\infty} (s^k_{\al ij}\psib_k\psi^j + s_{\al i})_L\psi^i_R \d_p^{\al}
$$
where $s^k_{\al i}, s^k_{\al ij}$ and $s_{\al i}$ are scalar functions of $x$. As in the proof of Proposition \ref{pJLO1} we look at the pseudodifferential order of these operators. The leading part $\psi^i_R\frac{\d}{\d x^i}$ of $\nabla$ has order $+1$, the two sub-leading terms have order $\leq 0$, while the remainder $s$ has order $\leq -1$. We calculate, for any $a\in \CS^0(M)$,
$$
[\i\eps\nabla, \rho(a)] = [\i\eps\nabla, (a\Pi)_L] = \i\eps\Big(\frac{\d a}{\d x^i}\Pi\Big)_{\!\! L}\psi^i_R  + \i\eps\Big(\Gamma^k_{ij}(x) p_k\frac{\d a}{\d p_j}\Pi\Big)_{\!\! L} \psi^i_R + \ldots
$$
We only write the terms of order $0$, and ignore the dots of order $-1$. In the same way 
$$
\nablab= \psib_{iR} \frac{\d}{\d p_i} + r
$$
has a leading term of order $-1$, and the remainder $r$ of order $-2$ can be expanded as $\sum_{|\al|=2}^{\infty} (r^i_{\al})_L\psib_{iR} \d_p^{\al}$ for some scalar functions $r^i_{\al}$. Hence 
$$
[\nablab,\rho(a)] = [\nablab, (a\Pi)_L] =\Big(\frac{\d a}{\d p_i}\Pi\Big)_{\!\! L}\psib_{iR} + \ldots
$$
is of order $-1$ and we ignore the dots of order $-2$. On the other hand, the generalized Laplacian $-D^2$ is given by (\ref{lap2}). Keeping only the leading terms we write
\beq
-D^2 &=& \i\eps \Big(\frac{\d}{\d x^i}\frac{\d}{\d p_i} +  (\Gamma^k_{ij})_L(\psi^i\psib_k)_R \frac{\d}{\d p_j}\Big) \non\\
&& + \frac{\eps^2 }{2}\, (\psi^i\psi^j)_R (R^k_{lij})_L  \Big( p_{kL} \frac{\d}{\d p_l} + (\psib_k\psi^l)_L\Big) + \ldots  \ , \non
\eeq
where the dots have the form of the leading terms but involve higher powers of the partial derivative $\d_p$ (hence have strictly lower order).  We proceed as in the proof of Proposition \ref{pJLO1} and consider $-D^2=\Delta + u$ as a perturbation of the flat Laplacian $\Delta = \i\eps \frac{\d}{\d x^i}\frac{\d}{\d p_i}$. A Duhamel expansion of the exponentials $\exp(-t_iD^2)$ appearing in the cochain $\varphi^D$ leads to the computation of terms like
\beq
\lefteqn{\Tr_s\big( \rho(a_0) \exp(t_0\Delta) \, X_1 \, \exp(t_1\Delta)\ldots X_k \, \exp(t_k\Delta) \big)=} \non\\
&& \qquad\qquad \bint \langle\!\langle \rho(a_0) \si^{t_0}_{\Delta}(X_1) \ldots \si^{t_0+\ldots +t_{k-1}}_{\Delta}(X_k) \, \exp \Delta  \rangle\!\rangle [n] \non
\eeq
where $X_i=u$ or $X_i=[D,\rho(a_j)]$ for some $a_j\in \CS^0(M)$. In particular $X_i$ has pseudodifferential order $\leq 0$, and this order is not modified by the action of the modular group $\si_{\Delta}$ because
$$
[\Delta, X_i] = \i\eps \Big( \frac{\d X_i}{\d x^j}\frac{\d}{\d p_j} + \frac{\d X_i}{\d p_j}\frac{\d}{\d x^j} + \frac{\d^2 X_i}{\d x^j\d p_j} \Big)\ .
$$ 
Now observe that in the above expressions for $-D^2$ and $[D,\rho(a)]$, a factor $\eps\psi_R$ always appears together with a pseudodifferential order $\leq 0$, whereas a factor $\psib_R$ always appears together with a pseudodifferential order $\leq -1$. The contraction map on the odd variables $\psi_R,\psib_R$ selects the only part of $\si^{t_0}_{\Delta}(X_1) \ldots \si^{t_0+\ldots +t_{k-1}}_{\Delta}(X_k)$ containing the product $(\psi^1\ldots \psi^n\psib_n\ldots \psib_1)_R$. This part has order $\leq -n$. Moreover, the dots in the above expressions for $-D^2$ and $[D,\rho(a)]$ contribute to an order $<-n$. A crucial consequence is that we only need to keep the leading terms of all quantities and ignore the dots because the Wodzicki residue vanishes on symbols of order $<-n$ (recall that the contractions $\langle\d_x^{\al}\d_p^{\beta}\exp\Delta\rangle$ do not affect the pseudodifferential order). Another crucial consequence is that all the derivatives $\d X_i/\d x^j$ appearing in the action of the modular group can be neglected, because these terms also contribute to an overall order $<-n$. Hence \emph{all functions of the variable $x$ behave like constants}. This drastically simplifies the computation of $\varphi^D$. One has
\beq
\lefteqn{\varphi_k^D(a_0,\ldots,a_k)= }\non\\
&& \int_{\Delta_k} \Tr_s \big(\rho(a_0)\si^{t_0}_{-D^2}([D,\rho(a_1)]) \ldots \si^{t_0+\ldots+t_{k-1}}_{-D^2}([D,\rho(a_k)])\exp(-D^2)\big) dt \non
\eeq
If we localize the supports of the symbols $a_i$ around a point $x_0\in U$ and choose a coordinate system in which $\Gamma^k_{ij}(x_0)=0$, we can write
\beq
[D,\rho(a)] &\simeq& \i\eps\Big(\frac{\d a}{\d x^i}\Pi\Big)_{\!\! L}\psi^i_R + \Big(\frac{\d a}{\d p_i}\Pi\Big)_{\!\! L}\psib_{iR} \ ,\non\\
-D^2 &\simeq& \Delta +  \frac{\eps^2 }{2} (\psi^i\psi^j)_R (R^k_{lij})_L  \Big( p_{kL} \frac{\d}{\d p_l} + (\psib_k\psi^l)_L\Big)  \non
\eeq
because we only keep the leading terms and ignore the $x$-derivatives of $\Gamma^k_{ij}$, hence $\Gamma^k_{ij}\simeq \Gamma^k_{ij}(x_0)=0$. For notational simplicity set $R^k_l = \frac{\eps^2}{2}(R^k_{lij})_L (\psi^i\psi^j)_R$ and recall that it behaves like a constant with respect to $x$. The generator of the  modular group $\si_{-D^2}$ is the commutator with $-D^2$. Its iterated actions on $X=[D,\rho(a)]$ read
\beq
&&-[D^2, X] \simeq [\Delta+ R^k_l  p_{kL} \frac{\d}{\d p_l} , X] \simeq \frac{\d X}{\d p_i} \Big(\i\eps\frac{\d}{\d x^i} + R^k_i p_{kL} \Big) \non\\
&&{[D^2,[D^2,X]]} \simeq \frac{\d^2 X}{\d p_i\d p_j} \Big(\i\eps\frac{\d}{\d x^i} + R^k_i p_{kL} \Big)\Big(\i\eps\frac{\d}{\d x^j} + R^l_j p_{lL} \Big) \non\\
&& \qquad \qquad \qquad +  R^j_i\frac{\d X}{\d p_i} \Big(\i\eps\frac{\d}{\d x^j} + R^l_j p_{lL} \Big)  \non
\eeq
Observe that the term $(R^k_{lij}\psib_k\psi^l)_L$ multiplied by $\rho(a)=(a\Pi)_L$ vanishes, because $R^k_{lij}\psib_k\psi^l \Pi = R^k_{lij}(\delta^l_k - \psi^l\psib_k) \Pi = R^k_{kij} \Pi$, and since $\Gamma$ is by hypothesis a Riemannian connection, $R^k_{kij}=0$. More generally
\beq
\si_{-D^2}^t(X) &=& \sum_{k=0}^{\infty} \frac{(-t)^k}{k!} \underbrace{[D^2,\ldots [D^2, }_k X ]\ldots ]\non\\
&\simeq& X + \sum_{k=1}^{\infty}\frac{t^k}{k!} \sum_{|\al|=1}^k P_{\al}(X) (\i \eps\d_x + p_L\cdot R)^{\al}  \non
\eeq
where $\al$ is a multi-index and $P_{\al}(X)$ is a linear combination of the partial $p$-derivatives of $X$. Since we drop the $x$-derivatives, the operator $(\i \eps\d_x + p_L\cdot R)^{\al}$ commutes with all operators under the graded trace so it can be moved to the right in front of $\exp(-D^2)$. Moreover $\rho(a)$ brings a factor $\Pi_L$ and we know that $R^k_{lij}\psib_k\psi^l\Pi=0$, so we may replace everywhere $-D^2$ by $\Delta + p_L\cdot R \cdot \d_p$. Then identities (\ref{iden}) lead to
\beq
\lefteqn{ \Tr_s\big( \rho(a_0) \si^{t_0}_{-D^2}([D,\rho(a_1)]) \ldots \si^{t_0+\ldots+t_{k-1}}_{-D^2}([D,\rho(a_k)]) \exp (\Delta + p_L\cdot R \cdot \d_p )\big)} \non\\
&& \qquad \qquad = \bint \big\langle\!\big\langle \rho(a_0) [D,\rho(a_1)] \ldots [D,\rho(a_k)] \exp (\Delta + p_L\cdot R \cdot \d_p ) \big\rangle\!\big\rangle[n]  \non
\eeq
The integral over $(t_0,\ldots,t_k)\in \Delta_k$ simply brings a factor $1/k!$. Lemma \ref{luse} applied to the matrix then gives
$$
\varphi_k^D(a_0,\ldots,a_k)= \frac{1}{k!} \bint \big\langle\!\big\langle \rho(a_0) [D,\rho(a_1)] \ldots [D,\rho(a_k)] \Td(R) \exp\Delta \big\rangle\!\big\rangle [n]\ .
$$
We have to select the coefficient of $\eps^n$ in this expression. $\eps$ always comes with a factor $\psi_R$ and the graded trace on the Clifford algebra selects the only polynomial $(\psi^1\ldots\psi^n\psib_n\ldots\psib_1)_R$, hence the variables $\psi_R$ and $\psib_{R}$ behave as if they anticommute. We make the identification with differential forms $\eps\psi^i_R\leftrightarrow dx^i$ and $\psib_{iR}\leftrightarrow dp_i-\Gamma^k_{ij}p_kdx^j$ over $T^*U$, which is consistent with the action of a coordinate change. Locally in our coordinate system one has $\Gamma^k_{ij}\simeq 0$ so that 
$$
[D,\rho(a)] \leftrightarrow \Big( \i\frac{\d a}{\d x^i} dx^i + \frac{\d a}{\d p_i} dp_i \Big) \Pi\ ,\quad \frac{\eps^2}{2} R^k_{lij}(\psi^i\psi^j)_R \leftrightarrow \frac{1}{2} R^k_{lij} dx^i\wedge dx^j=R^k_l\ .
$$ 
To be more precise, if we multiply the bracket by the volume form of the cotangent bundle $\om^n/n!=dp_1\wedge dx^1 \ldots dp_n\wedge dx^n$, and compare it to the normalization condition $\langle(\psib_1\psi^1\ldots \psib_n\psi^n)_R\rangle =(-1)^n$, one finds the equality of $2n$-forms over $T^*U$ (the subscript $_{\mathrm{vol}}$ denotes the top-component of a differential form) 
\beq
\lefteqn{\big\langle\!\big\langle \rho(a_0) [D,\rho(a_1)] \ldots [D,\rho(a_k)] \Td(R) \exp\Delta \big\rangle\!\big\rangle [n] \, \frac{\om^n}{n!} = }\non\\
&& \qquad\qquad (-1)^n \i^{k-n}\, \big( a_0 da_1  \ldots da_k \Td(R) \Pi \big)_{\mathrm{vol}}     +\ \mbox{terms of order}\ < -n \non
\eeq
The first term of the right-hand-side is a scalar symbol of order $\leq -n$, times the volume form. We claim that this symbol in fact has order $<-n$. Indeed the product $a_0 da_1  \ldots da_k$ brings $n$ partial derivatives with respect to the variables $( p_1,\ldots p_n)$. Writing its leading symbol in polar coordinates $(\|p\|, \theta_1,\ldots,\theta_{n-1})$, one sees that is is proportional to $\|p\|^{1-n}$ times a partial derivative $\frac{\d a}{\d \|p\|}$. The latter has order $\leq -2$. Hence the Wodzicki residue vanishes. \cqfd\\

We now deal with the Radul cocycle. Let $q\in \CS^1(M)$ be a symbol of order one, with positive and invertible leading symbol. The logarithm $\log q$ is no longer classical, but belongs to the larger class of log-polyhomogeneous symbols: its asymptotic expansion in a local coordinate system $(x,p)$ reads
\be
(\log q)(x,p)= \log\|p\| + q'_0(x,p)
\ee
where $q'_0\in\CS^0(M)$ is a classical symbol of order $\leq 0$. It is easy to check that the commutator (for the $\star$-product) of $\log q$ with any classical symbol $a\in \CS^m(M)$ is in $\CS^{m-1}(M)$. In fact $[\log q, a]$ has an expansion
\be
[\log q, a] = \sum_{k=1}^{\infty} \frac{1}{k}(-1)^{k-1} a^{(k)} q^{-k} \label{log}
\ee
where $a^{(k)}\in \CS^m(M)$ denotes the $k$-th power of the derivation $[q,\ ]$ on $a$. Thus $[\log q,\ ]$ is an outer derivation on the algebra of classical symbols $\CS(M)$. The Radul cocycle \cite{Ra} is the bilinear map $c:\CS(M)\times\CS(M)\to\cc$ defined by means of the Wodzicki residue 
\be
c(a_0,a_1) = \bint a_0  [\log q,a_1 ]\ ,\quad \forall a_i\in \CS(M)\ .
\ee
The expansion (\ref{log}) shows that the Wodzicki residue vanishes on commutators $[\log q,a]$ for any classical symbol $a$. Hence the Wodzicki residue is trace on $\CS(M)$ which is closed with respect to the derivation $[\log q,\ ]$. Elementary algebraic manipulations show the antisymetry property $c(a_0,a_1)=-c(a_1,a_0)$. Moreover the Hochschild coboundary of $c$ is
$$
bc (a_0,a_1,a_2)= c(a_0a_1, a_2) - c(a_0,a_1a_2) + c(a_2a_0,a_1)=0
$$
for all $a_i\in \CS(M)$. Thus $c$ is a cyclic one-cocycle. Originally $c$ was introduced as a two-cocycle over the Lie algebra $\CS(M)$, with commutator as Lie bracket, but the cyclic cocycle property is actually stronger. From now on we view $c$ as a cyclic one-cocycle over the subalgebra $\CS^0(M)\subset\CS(M)$ of symbols of order $\leq 0$.\\
Then we extend the commutator $[\log q,\ ]$ to a derivation on the algebra $\Lc(M)\subset \End(\CS(M,E))$ as follows. Recall that $\Lc(M)$ is generated by left multiplications $a_L$ for all symbols $a\in \CS(M,E)$, and right multiplications $b_R$ for all polynomial symbols $b\in \PS(M,E)$. Then extend $q\in \CS^1(M)$ to an elliptic positive symbol $\tilde{q}\in \CS^1(M,E)$ of scalar type and set
\be
\delta(a_Lb_R)=([\log \tilde{q},a])_Lb_R \qquad \forall\  a \in \CS(M,E)\ ,\ b\in \PS(M,E)\ .
\ee
Since the left representation $a\mapsto a_L$ is faithful, $\delta:\Lc(M)\to\Lc(M)$ is well-defined. It is clearly a derivation. In an obvious fashion we extend it to a derivation, still denoted $\delta$, on the algebra of formal power series $\Sc(M)=\Lc(M)[[\eps]]$ by setting $\delta\eps=0$. It has good properties with respect to the subspaces $\Dc^m_k(M)$. Indeed in a local coordinate system over $U\subset M$, one has
\beq
\delta(\d_p) &=& \i \delta(x_R-x_L)\ =\ -\i ([\log \tilde{q},x])_L\ \in \CS^{-1}(U,E)_L\non\\
\delta(\d_x) &=& \i \delta(p_L-p_R)\ =\ \i ([\log \tilde{q},p])_L\ \in \CS^{0}(U,E)_L
\eeq
and also $\delta(\CS^{m}(M,E)_L)\subset \CS^{m-1}(M,E)_L$ for all $m\in\rr$. This shows that $\delta(\Dc^m_k(M))\subset \Dc^{m-1/2}_k(M)$ for all $m\in\rr$ and $k\in\nn$. If $\Delta\in \Dc^{1/2}_1(M)$ is a generalized Laplacian, one has
$$
\delta \exp(\Delta) = \int_0^1 e^{t\Delta}\, \delta\Delta\, e^{(1-t)\Delta}\, dt = \int_0^1 \si_{\Delta}^t(\delta\Delta) \exp (\Delta) \, dt
$$
hence $\delta$ restricts to a derivation on the $\Dc(M)$-bimodule $\Tc(M)$ of trace-class operators. The analogue of expansion (\ref{log}) for $\delta$ shows that the graded trace $\Tr_s:\Tc(M)\to\cc$ is $\delta$-closed.

\begin{proposition}
Let $D\in \Dc^1(M)$ be a generalized Dirac operator and $\delta$ the derivation associated to an elliptic positive symbol $\tilde{q}\in \CS^1(M,E)$. The homogeneous cochains over the algebra $\CS^0(M)$
\beq
\lefteqn{\varphi_k^{D,\delta}(a_0,\ldots,a_k)= }\\
&& \hspace{-0.5cm} \sum_{i=1}^k(-1)^{i+1}\int_{\Delta_k} \!\!\! \Tr_s \big( \rho(a_0)e^{-t_0D^2}[D,\rho(a_1)] e^{-t_1D^2} \ldots \delta\rho(a_i)e^{-t_iD^2} \ldots [D,\rho(a_k)]e^{-t_kD^2}\big) dt \non\\
&& \hspace{-0.7cm} +\sum_{i=1}^{k+1} (-1)^{i} \int_{\Delta_{k+1}}\!\!\! \Tr_s \big(\rho(a_0)e^{-t_0D^2}[D,\rho(a_1)] e^{-t_1D^2} \ldots \delta D e^{-t_iD^2} \ldots [D,\rho(a_k)]e^{-t_kD^2}\big) dt \non
\eeq
defined for all $k\in 2\nn+1$, are the components of an odd periodic cyclic cocycle $\varphi^{D,\delta}$ and vanish whenever $k>2n+1$, $n=\dim M$. Moreover, the periodic cyclic cohomology class $[\varphi^{D,\delta}]\in HP^1(\CS^0(M))$ does not depend on $D$ nor $\tilde{q}$. 
\end{proposition}
{\it Proof:} Analogous to Proposition \ref{pJLO1}. Details are left to the reader. \cqfd\\

\begin{proposition}
Let $D=-\i\eps d_R + \nablab$ be a de Rham-Dirac operator. Then the first component $\varphi^{D,\delta}_1=c$ is the Radul cocycle on $\CS^0(M)$, while the other components $\varphi^{D,\delta}_k$ vanish for $k>1$. Hence $[\varphi^{D,\delta}]$ is the periodic cyclic cohomology class of $[c]$. 
\end{proposition}
{\it Proof:} We proceed as in Proposition \ref{pwod}. The commutator $[D,\rho(a)]$ only brings $\psib_R$ which cannot be balanced by $\psi_R$, hence $\varphi_k^{D,\delta}$ vanishes whenever $k>1$. The only non-zero component is
$$
\varphi_1^{D,\delta}(a_0,a_1) = \int_0^1 \Tr_s \big(\rho(a_0) e^{-t D^2} \delta\rho(a_1) e^{(t-1)D^2} \big) dt\ .
$$
Observe that
$$
\frac{d}{dt}\, \Tr_s \big(\rho(a_0) e^{-t D^2} \delta\rho(a_1) e^{(t-1)D^2} \big) = -\Tr_s \big(\rho(a_0) e^{-t D^2} [D^2,\delta\rho(a_1)] e^{(t-1)D^2} \big)\ .
$$
The identity $[D^2,\delta\rho(a_1)]=D[D,\delta\rho(a_1)]+ [D,\delta\rho(a_1)]D$ and the graded trace property yield
$$
-\Tr_s \big(\rho(a_0) e^{-t D^2} [D^2,\delta\rho(a_1)] e^{(t-1)D^2} \big) = \Tr_s \big([D,\rho(a_0)] e^{-t D^2} [D,\delta\rho(a_1)] e^{(t-1)D^2} \big)
$$
This quantity vanishes because the commutators $[D,\rho(a)]$ are proportional to $\psib_R$. Hence $\Tr_s \big(\rho(a_0) e^{-t D^2} \delta\rho(a_1) e^{(t-1)D^2} \big)$ does not depend on $t$ and we can rewrite the integral $\varphi_1^{D,\delta}$ in terms of its integrand at $t=0$:
$$
\varphi_1^{D,\delta}(a_0,a_1) = \Tr_s \big(\rho(a_0)\delta\rho(a_1) e^{-D^2} \big) = \Tr_s \big((a_0 [\log q, a_1]\Pi)_L e^{-D^2} \big)\ .
$$
The computation is now completely analogous to Proposition \ref{pwod} and one finds
$$
\varphi_1^{D,\delta}(a_0,a_1) = \bint a_0 [\log q, a_1]
$$
as claimed. \cqfd\\

Choose an affine torsion-free connection on the tangent bundle $TM$, and let $R\in \Om^2(M,\End(TM))$ be its curvature two-form. The Todd class of the complexified tangent bundle $\Td(T_{\cc}M)\in H^{\bullet}(M,\cc)$ is the cohomology class of even degree represented by the closed differential form 
\be
\Td(\i R/2\pi) = \det\left(\frac{\i R/2\pi}{e^{\i R/2\pi}-1}\right)  \ ,
\ee
where the determinant acts on the sections of the endomorphism bundle of $T_{\cc}M$.

\begin{theorem}\label{trad}
Let $M$ be a closed manifold. The periodic cyclic cohomology class of $[c]\in HP^1(\CS^0(M))$ is
\be
[c] =  \la^*\big([S^*M]\cap\, \pi^*\Td(T_{\cc}M)\big)\ ,
\ee
where $\la^*$ is the pullback (\ref{iso}) induced by the leading symbol homomorphism, $\Td(T_{\cc}M) \in H^{\bullet}(M,\cc)$ is the Todd class of the complexified tangent bundle,  and $\pi:S^*M\to M$ is the cosphere bundle endowed with its canonical orientation and fundamental class $[S^*M]\in H_{\bullet}(S^*M)$. 
\end{theorem}
{\it Proof:} We apply verbatim the proof of Theorem \ref{twod}. We can replace the commutator $[D,\rho(a)]$ by $\i\eps (\frac{\d a}{\d x^i}\Pi )_L \psi^i_R + (\frac{\d a}{\d p_i}\Pi)_L \psib_{iR} $ in a local coordinate system, and consider $R^k_l = \frac{\eps^2}{2}(R^k_{lij})_L (\psi^i\psi^j)_R$ as independent of $x$. Then
\beq
\lefteqn{\varphi_k^{D,\delta}(a_0,\ldots,a_k)= } \non\\
&&  \sum_{i=1}^k \frac{(-1)^{i+1}}{k!} \bint \big\langle\!\big\langle \rho(a_0) [D,\rho(a_1)] \ldots \delta\rho(a_i) \ldots [D,\rho(a_k)] \Td(R)\exp\Delta \big\rangle\!\big\rangle [n]  \non \\
&& + \sum_{i=1}^{k+1} \frac{(-1)^i}{(k+1)!} \bint \big\langle\!\big\langle \rho(a_0) [D,\rho(a_1)] \ldots \delta D \ldots [D,\rho(a_k)] \Td(R)\exp\Delta \big\rangle\!\big\rangle [n]  \non
\eeq
The first term of the right-hand-side vanishes. Indeed, the bracket selects the polynomial $(\psib_1\ldots\psib_n)_R$ which brings $n$ derivatives with respect to $p$, and $\delta\rho(a)=([\log q,a]\Pi)_L$ is of order $-1$. Hence the symbol under the Wodzicki residue has order $< -n$ and disappears. We are left with the second term involving $\delta D$. Recall that $(\log q) (x,p)= \log \|p\| + q'_0(x,p)$ where $q'_0$ is a classical symbol of order $\leq 0$. At leading order one has
\beq
\delta D &=& -\i\eps\Big( \frac{\d \log q}{\d x^i}\Big)_{\!\! L} \psi_R^i-\Big( \frac{\d \log q}{\d p_i}\Big)_{\!\! L} \psib_{iR} + \ldots \non\\
&=&  -\i\eps\Big( \frac{\d q'_0}{\d x^i}\Big)_{\!\! L} \psi_R^i - \Big(\frac{\d q'_0}{\d p_i} + \frac{p^i}{\|p\|^2} \Big)_{\!\! L} \psib_{iR} + \ldots \non
\eeq
where $p^i =  \delta^{ij}p_j$. The leading term proportional to $\psi_R$ (resp. $\psib_R$) is of order $\leq 0$ (resp. $\leq -1$), and the dots proportional to $\psi_R$ (resp. $\psib_R$) are of order $<0$ (resp. $<-1$). The bracket under the residue is expressed by means of differential forms:
\beq
\lefteqn{\big\langle\!\big\langle \rho(a_0) [D,\rho(a_1)] \ldots \delta D \ldots [D,\rho(a_k)] \Td(R)\exp\Delta \big\rangle\!\big\rangle [n]\, \frac{\om^n}{n!}= }\non\\
&&  -(-1)^n\i^{k+1-n} \Big( a_0 da_1  \ldots \Big(dq'_0 + \frac{p^idp_i}{\|p\|^2}\Big) \ldots da_k \Td(R) \Pi \Big)_{\mathrm{vol}}  +\ldots \non
\eeq
The leading part is a symbol of order $\leq -n$, while the dots of order $<-n$ are killed by the Wodzicki residue. One shows as in the proof of \ref{twod} that the term $a_0 da_1  \ldots dq'_0 \ldots da_k$ is also killed. Hence the only remaining term is proportional to $p^idp_i/\|p\|^2$. At leading order we can view $a_0,\ldots,a_k$ as scalar functions over the cosphere bundle. Since $\tr_s(\Pi)=1$ the residue becomes the integral of a $(2n-1)$-form (remark that it is globally defined)
\beq
\varphi_k^{D,\delta}(a_0,\ldots,a_k) &=& \frac{(-1)^n\i^{k+1-n}}{(2\pi)^n k!} \int_{S^*M} \iota(L)\cdot \Big(\frac{p^idp_i}{\|p\|^2}\wedge a_0 da_1 \ldots da_k \Td(R) \Big)\non \\
&=& \frac{\i^{k+n+1}}{(2\pi)^n k!} \int_{S^*M}  a_0 da_1 \ldots da_k \Td(R) \non 
\eeq
where $L=p_i\frac{\d}{\d p_i}$ is the fundamental vector field on $T^*M$. The dimension of $S^*M$ equals $2n-1$ and the parity of the cochain is actually odd, so one gets
$$
\varphi_{2k+1}^{D,\delta}(a_0,\ldots,a_{2k+1}) = \frac{1}{(2\pi\i)^{k+1} (2k+1)!} \int_{S^*M}  a_0 da_1 \ldots da_{2k+1} \Td(\i R/2\pi) 
$$
for any $k\in\nn$. This is precisely the pullback, under the morphism $\la$, of the degree $2k+1$ component of the de Rham cycle $[S^*M]\cap \Td(\i R/2\pi)$. \cqfd\\

\section{Atiyah-Singer index theorem}\label{sit}

An immediate corollary of Theorem \ref{trad} is the Atiyah-Singer index theorem, which computes the index of an elliptic pseudodifferential operator on a closed manifold $M$, in terms of local data. We consider the algebra $\CL^0(M)$ of scalar pseudodifferential operators of order $\leq 0$ as an extension of the algebra $\CS^0(M)$ of formal symbols, with kernel the algebra of smoothing operators: 
\be
(E):\quad 0 \to \L^{-\infty}(M) \to \CL^0(M) \to \CS^0(M) \to 0\ .
\ee
An operator $Q\in \CL^0(M)$ is elliptic if and only if its leading symbol is invertible, or equivalently, if its formal symbol is invertible in $\CS^0(M)$. Thus $Q$ has a parametrix $P\in \CL^0(M)$ which is an inverse modulo smoothing operators, that is $PQ-1$ and $QP-1$ are in $\L^{-\infty}(M)$. The obstruction of perturbing $Q$ to an exactly invertible operator in $\CL^0(M)$ is measured by the \emph{index map} of the extension $(E)$ in algebraic $K$-theory
\be
\Ind_E: K_1(\CS^0(M)) \to K_0(\L^{-\infty}(M))\cong\zz \label{ind}
\ee
cf. \cite{M}. Indeed the formal symbol of $Q$ is invertible hence defines a class $[Q]$ in the algebraic $K$-theory group $K_1(\CS^0(M))$, and its image under the map (\ref{ind}) coincides with the Fredholm index of $Q$ as a bounded operator on $L^2(M)$. In \cite{P8} we presented a general procedure allowing to compute local index formulas associated to extensions. In the simple case of pseudodifferential operators on a closed manifold, the calculation of the index reduces to the Radul cocycle evaluated on $Q$ and its parametrix (more precisely, on their formal symbols):
\be
\Ind_E([Q]) = c(P,Q)\ .
\ee
In terms of Connes' pairing between $K$-theory and cyclic cohomology \cite{C86}, the above formula is precisely the pairing of $[Q]\in K_1(\CS^0(M))$ with the cyclic cohomology class $[c]\in HP^1(\CS^0(M))$. Since by Theorem \ref{trad}, this class is a pullback under the leading symbol map $\la:\CS^0(M)\to\cinf(S^*M)$, we are able to express the index of $Q$ in terms of its leading symbol which is an invertible function $g\in\cinf(S^*M)$. This is not surprising because the algebra $\CS^0(M)$ is a pro-nilpotent extension of $\cinf(S^*M)$, and this implies an isomorphism of the algebraic $K$-theory groups $K_1(\CS^0(M))\cong K_1(\cinf(S^*M))$. In fact one has a diagram of extensions
$$
\xymatrix{
0 \ar [r] & \L^{-\infty}(M) \ar[r] \ar[d] & \CL^0(M) \ar[r] \ar@{=}[d] & \CS^0(M)  \ar[r] \ar[d] & 0  \\
0 \ar [r] & \CL^{-1}(M) \ar[r] & \CL^0(M) \ar[r] & \cinf(S^*M)  \ar[r] & 0  }
$$
The vertical arrows are isomorphisms both at the $K$-theoretic and periodic cyclic cohomology levels. Thus the index map of $(E)$ should really be viewed as a map 
\be
\Ind_E: K_1(\cinf(S^*M))\to \zz\ ,
\ee
sending the leading symbol class $[g]\in K_1(\cinf(S^*M))$ to the Fredholm index of $Q$. Of course, everything extends to pseudodifferential operators acting on the sections of a (trivially graded) complex vector bundle over $M$, the leading symbols being matrix-valued functions over $S^*M$. In order to state the index formula we need to recall that any class $[g]\in K_1(\cinf(S^*M))$, represented by an invertible matrix-valued function $g$, has a Chern character in the cohomology $H^{\bullet}(S^*M,\cc)$ of odd degree represented by the closed differential form
\be
\ch(g) = \sum_{k\geq 0}  \frac{k!}{(2k+1)!} \, \tr\Big(\frac{(g^{-1}dg)^{2k+1}}{(2\pi\i)^{k+1}}\Big)\ .
\ee
\begin{corollary}[Index theorem] 
Let $Q$ be an elliptic pseudodifferential operator of order $\leq 0$ acting on the sections of a trivially graded vector bundle over $M$, with leading symbol class $[g]\in K_1(\cinf(S^*M))$. Then the Fredholm index of $Q$ is the integer
\be
\Ind(Q) = \langle [S^*M], \pi^*\Td(T_{\cc}M)\cup \ch([g])\rangle\ .
\ee
\end{corollary}
{\it Proof:} If $\varphi=(\varphi_1,\varphi_3,\ldots,\varphi_{2n-1})$ is an odd $(b+B)$-cocycle over $\CS^0(M)$, its pairing with the $K$-theory class $[Q]\in K_1(\CS^0(M))$ reads (\cite{C94})
$$
\langle [\varphi] , [Q] \rangle = \sum_{k\geq 0} (-1)^k\,k!\, (\varphi_{2k+1}\otimes \tr) (P,Q,\ldots,P,Q)
$$
where, strictly speaking, $Q$ and its parametrix $P$ should be replaced by their formal symbols. If $\varphi$ is the pullback of an odd homology class $[C]\in H_{\bullet}(S^*M,\cc)$ under the leading symbol map $\la$, the above formula factors through the leading symbols $g=\la(Q)$ and $g^{-1}=\la(P)$. Using the identity $d g^{-1}=-g^{-1}dg g^{-1}$ one gets
$$
\langle [\varphi] , [Q] \rangle = \sum_{k\geq 0} \frac{(-1)^k k!}{(2\pi\i)^{k+1}(2k+1)!} \, \langle C_{2k+1}, \tr(g^{-1}dg (dg^{-1}dg)^k) \rangle = \langle [C] , \ch([g]) \rangle\ .
$$
Applying this formula to the periodic cyclic cohomology class of $c$ given by Theorem \ref{trad} gives the desired formula for $\Ind(Q) = \langle [c] , [Q] \rangle$. \cqfd \\

\end{document}